\documentclass[%
  twocolumn
   , colorlinks % also like hyperref with CSC colorscheme (this is the default)
]{mpi2015-cscpreprint}

\usepackage[american]{babel}
\usepackage[protrusion=true, expansion=true]{microtype}

\usepackage{graphicx}

\usepackage{amssymb}
\usepackage{amsthm}

\usepackage{wrapfig}
\usepackage{tikz}
\usepackage{pgfplots}
\pgfplotsset{compat=1.18}
\usepackage{filecontents}
\usepackage{subcaption}
\usepackage{caption}
\usepackage{wrapfig}
\usepackage{cuted}
\usepackage{hyperref}
\usepackage[linesnumbered,ruled,vlined]{algorithm2e}
\newcommand{\norm}[1]{\left\lVert#1\right\rVert}

\usepackage{xcolor}
\definecolor{oiBlue}{HTML}{0072B2}
\definecolor{oiOrange}{HTML}{E69F00}
\definecolor{oiGreen}{HTML}{009E73}
\usepackage[export]{adjustbox} % in preamble
\usepackage{cuted}
\AtBeginDocument{

}

\newcommand{\blankgraphics}[2]{% #1 = width, #2 = height
  \begingroup
  \setlength{\fboxsep}{0pt}%
  \fbox{\rule{#1}{0pt}\rule{0pt}{#2}}% frame + empty box
  \endgroup
}

\begin{document}

%%%%%%%%%%%%%%%%%%%%%%%%%%%%%%%%%%%%%%%%%%%%%%%%%%%%%%%%%%%%%%%%%%%%%%%%%%%%%%%%
% PAPER INFORMATION.                                                           %
%%%%%%%%%%%%%%%%%%%%%%%%%%%%%%%%%%%%%%%%%%%%%%%%%%%%%%%%%%%%%%%%%%%%%%%%%%%%%%%%

\title{On the Compatibility of DEIM Subspaces for Nonlinear Forces Reduction in Projective Dynamics}

\author[$\ast$]{Shaimaa Monem$^\dagger$}
\affil[$\ast$]{Max Planck Institute for Dynamics of Complex Technical Systems, Sandtorsraße 1, 19306 Magdeburg, Germany.\authorcr%
  \email{monem@mpi-magdeburg.mpg.de}, \orcid{0009-0008-4038-3452}}
  \affil[$\dagger$]{Corresponding author. Permanent email:\authorcr%
  \email{monem8shaimaa@gmail.com}}% chktex 8

\author[$\ast$]{Igor Pontes Duff}
\affil[$\ast$]{Max Planck Institute for Dynamics of Complex Technical Systems, Sandtorsraße 1, 19306 Magdeburg, Germany.\authorcr%
\email{pontes@mpi-magdeburg.mpg.de}, \orcid{0000-0001-6433-6142}}% chktex 8

\author[$\ast \ddagger$]{Peter Benner}
\affil[$\ast$]{Max Planck Institute for Dynamics of Complex Technical Systems, Sandtorsraße 1, 19306 Magdeburg, Germany.\authorcr%
\email{benner@mpi-magdeburg.mpg.de}, \orcid{0000-0003-3362-4103}}% chktex 8
\affil[$\ddagger$]{Otto-von-Guericke Universität, Magdeburg, Germany.}

\shorttitle{DEIM for PD}
\shortauthor{Monem et. al.}
\shortdate{}



\abstract{
	In physics-based simulations, the evaluation of nonlinear constraint forces represents a primary source of computational complexity. Within the framework of projective dynamics, the solver alternates between estimating these forces and updating vertex positions for the subsequent frame. This work investigates the potential of snapshot-based subspace methods to accelerate the computation of constraint terms while maintaining the fidelity of nonlinear dynamics. We focus specifically on the Discrete Empirical Interpolation Method (DEIM), examining whether it can identify a reduced set of constrained mesh elements that enable accurate interpolation and approximation of forces across the full system.

	Our experiments demonstrate that DEIM is neither intuitive nor consistently effective for constraint projection reduction. In particular, its performance degrades when constraints involve multiple vertices. To validate these observations, we present comparisons between full-order and reduced-order simulations on both surface meshes embedded in three dimensions and volumetric meshes. Furthermore, we report reconstruction errors for various classes of constraints, providing quantitative evidence of the method’s limitations in this context.
}

\novelty{Recent research has demonstrated that snapshot-based reduced subspaces for vertex positions outperform traditional skinning-based approaches in both accuracy and computational efficiency. Unlike skinning methods, snapshot-based techniques can capture a wider range of motions and deformations while more accurately preserving rotational effects. In this context, the Discrete Empirical Interpolation Method (DEIM) emerges as a promising candidate for accelerating the evaluation of nonlinear forces. This raises a central research question: to what extent can DEIM alleviate the computational bottleneck associated with nonlinear terms while preserving their essential physical behavior?
}

\maketitle

%%%%%%%%%%%%%%%%%%%%%%%%%%%%%%%%%%%%%%%%%%%%%%%%%%%%%%%%%%%%%%%%%%%%%%%%%%%%%%%%
% PAPER CONTENT.                                                               %
%%%%%%%%%%%%%%%%%%%%%%%%%%%%%%%%%%%%%%%%%%%%%%%%%%%%%%%%%%%%%%%%%%%%%%%%%%%%%%%%

\section{Introduction}%
\label{sec:intro}

	The simulation of deformable objects plays a fundamental role in computer graphics and related fields, enabling the realistic representation of materials such as cloth, soft tissues, elastic solids, and thin shells. The plausibility of deformable models is crucial for applications ranging from visual effects in film to interactive experiences in video games, where physical consistency enhances user engagement. Beyond visual fidelity, these simulations provide animators and designers with physically grounded tools for manipulating characters and environments, effectively bridging the gap between artistic intent and physical accuracy—particularly in real-time, interactive scenarios that require faithful representation of underlying forces.

	In scientific and engineering domains, deformable simulations are equally essential. Biomedical applications rely on soft tissue simulation for surgical planning, medical training, and prosthetics design, where nonlinear and anisotropic material behavior must be accurately captured \cite{Zhang19,PBSSurvay14}. Similarly, mechanical and civil engineering employ deformable models to analyze stress and strain in flexible components, safety systems, and load-bearing structures. The ability to reproduce deformation dynamics in real time has become a cornerstone of virtual prototyping, substantially reducing the need for costly physical testing. Recent advances in computational techniques have significantly broadened the scope of deformable simulations. Methods such as finite element analysis \cite{snowden2017methods}, position-based dynamics \cite{PBS15,PBSSurvay14}, and projective dynamics \cite{Bouaziz14} offer complementary trade-offs between robustness, computational speed, and physical fidelity. Projective Dynamics, in particular, has emerged as a powerful framework that reformulates the evolution of physical systems as a sequence of constraint projections, enabling stable and efficient simulation even for highly nonlinear deformations \cite{Bouaziz14}. To further accelerate these methods, model order reduction has also been explored, compressing high-dimensional nonlinear dynamics into compact, low-dimensional subspaces. These techniques have been successfully applied across engineering \cite{GKIMISIS2025118115}, fluid dynamics \cite{lassila14}, computer graphics \cite{Trusty23}, gaming \cite{Bnneland19}, and bio-systems \cite{Radulescu12}, highlighting the cross-disciplinary relevance of deformable simulation as both a scientific and creative tool.

	The core computational bottleneck in projective dynamics lies in constraint projections, which typically account for $70\%$ or more of the total computational time. This work presents the first attempt to integrate DEIM as a subspace-based approach for approximating nonlinear constraint forces within the projective dynamics framework.

%%%%%%%%%%%%%%%%%%%%%%%%%%%%%%%%%%%%%%%%%%%%%%%%%%%%%%%%%%%%%%%%%%%%%%%%%%%%%%%%
%%%%%%%%%%%%%%%%%%%%%%%%%%%%%%%%%%%%%%%%%%%%%%%%%%%%%%%%%%%%%%%%%%%%%%%%%%%%%%%%

\section{Related Work}
\label{related_work}

	%%%%%%%%%%%%%%%%%%%%%%%%%%%%%%%%%%%%%%%%%%%%%%%%%%%%%%%%%%%%%%%%%%%%%%%%%%%%%%%%
	\subsection{Projective Dynamics}

	Projective Dynamics offers a fast and stable solver for real-time simulation of cloth, soft bodies, shells, and example-based animations. It combines quadratic energy potentials with a local-global alternating optimization scheme and has become one of the most influential real-time simulation frameworks in graphics and animation.

	Projective Dynamics bridges two established approaches: Position-Based Dynamics (PBD) \cite{PBS15, DKWB2018, MMMCN16}, which performs constraint projections per timestep and is widely adopted in games and interactive applications due to its robustness; and the implicit Finite Element Method (FEM) with Newton–Raphson solvers \cite{snowden2017methods, BDWA98}, known for its high accuracy and robustness but hindered by the cost of matrix factorizations.

	Projective Dynamics supports a wide range of applications in computer graphics, making it a continual target for optimization and acceleration. Its structure can be interpreted as a quasi-Newton method \cite{LTB17}, allowing support for hyperelastic and spline-based materials. When combined with the state-based peridynamics framework, it can be extended to simulate elastoplastic materials \cite{Xiaowei18}, and has also been applied in differentiable soft-body simulators \cite{DiffPD21}. Incorporating adaptive smoothing within a full multigrid framework, along with GPU acceleration, has been shown to significantly improve convergence rates \cite{Zhendong18}.

	To further accelerate simulations, reduced subspaces have been introduced for both vertex position updates and constraint projections \cite{Brandt18, Brandt19, Monem23}. In this work, we investigate the possibility of selecting a small subset of mesh elements to approximate the nonlinear forces acting on the remaining constrained elements using the discrete empirical interpolation method. This study is motivated by the superior performance of snapshot-based subspaces \cite{Monem23, Monem25} compared to linear blend skinning subspaces \cite{Brandt18, Brandt19} in capturing vertex motion.

	%%%%%%%%%%%%%%%%%%%%%%%%%%%%%%%%%%%%%%%%%%%%%%%%%%%%%%%%%%%%%%%%%%%%%%%%%%%%%%%%
	\subsection{Discrete Empirical Interpolation Method}

	A widely used approach for accelerating the simulation of deformable elastic bodies is to construct a reduced subspace spanned by a limited set of basis vectors. This subspace contains significantly fewer degrees of freedom while still approximating the behavior of the original high-dimensional system. Computational efficiency is achieved by projecting the governing equations onto this reduced space, where all subsequent computations are performed \cite{huang2019survey}.

	Methods for identifying suitable bases have been extensively explored in engineering, computer graphics, and animation \cite{Pentland_Williams89, Barbivc11, Barbic05}. The foundational work of Terzopoulos et al.~\cite{Terzopoulos87} established a theoretical framework for elastic deformable materials, which later inspired reduced-order methods in graphics. Pentland and Williams \cite{Pentland_Williams89} introduced vibration-based reduction subspaces for animation, marking the introduction of modal analysis into computer graphics. These early modal bases were computed from mass and stiffness matrices, relying solely on the undeformed rest state \cite{Pentland_Williams89}.

	More recently, snapshot-based techniques have been shown to strike an effective balance between accuracy and computational cost for reducing vertex position computations \cite{Monem23}. In this study, we investigate DEIM as a candidate snapshot-based method for constructing reduced subspaces aimed at accelerating nonlinear constraint projections.

	Various extensions of the Discrete Empirical Interpolation Method (DEIM) have been developed, with numerous studies investigating sampling adaptivity, stability, and error bound refinement \cite{Chaturantabut10, Drmavc16, Benjamin20, chaturantabut2012state}. DEIM originated as the discrete counterpart of the Empirical Interpolation Method (EIM) \cite{Barrault04}, which focused on reduced bases for partial differential equations arising from finite element discretizations.

	The core idea of DEIM is to identify a small set of spatial points at which a generic nonlinear discretized function is interpolated. These points are selected based on where the subspace projection error is highest, and exact interpolation is enforced at those locations. This strategy reduces the approximation error and enables faster evaluation of the nonlinear components of the full model.

	DEIM has been effectively combined with Proper Orthogonal Decomposition (POD) to reduce the computational cost of large-scale models. This combination has been successfully applied across various domains, including nonlinear viscous fingering in porous media \cite{Chaturantabut10}, implicit shallow water models for geophysical flows \cite{cstefuanescu2013pod}, the dynamic simulation of two-dimensional catalytic reactor models \cite{bremer2017pod}, and was applied in application to generic parabolic PDEs \cite{maday2008parareal}. Beyond standard applications, DEIM has also been utilized in advanced tasks such as sensor placement optimization for nuclear reactors \cite{argaud2018sensor}.

	\begin{itemize}
		\item We conduct a systematic evaluation of DEIM across different types of constraint projections on simplified examples, assessing its suitability for accelerating physics-based animation.
		\item We present both qualitative (visual) and quantitative (error-based) comparisons between full-order and reduced-order simulations on surface and volumetric meshes, highlighting the strengths and limitations of DEIM in this context.
	\end{itemize}
	
	It is worth noting that recent advances in hyper-reduction have introduced novel frameworks, such as Deep-HyROMnet for nonlinear parameterized PDEs \cite{cicci2022deep} and Hyper-Reduced Autoencoders, which combine neural networks with collocation based hyper reduction techniques \cite{cocola2023hyper}. However, exploring the integration of these approaches into projective dynamics lies beyond the scope of this work.
	%%%%%%%%%%%%%%%%%%%%%%%%%%%%%%%%%%%%%%%%%%%%%%%%%%%%%%%%%%%%%%%%%%%%%%%%%%%%%%%%

	\section{Background}

	The motion of vertex positions in an elastic mesh is governed by Newton's second law. In the absence of damping forces, and under the influence of stiffness, internal, and external forces—with corresponding matrices $\mathbf{K}$ and $\mathbf{F}$—the elastic model is expressed as

	$$
	\mathbf{M} \ddot{q} + \mathbf{K} q = \mathbf{F},
	$$ where $\mathbf{M}$ is the lumped mass matrix, and $q(t) \in \mathbb{R}^{N \times 3}$ denotes the vertex position displacements as a function of time.

	\begin{figure*}[t!]
		\centering
		\includegraphics[width=\textwidth, height=0.4\textwidth, keepaspectratio]{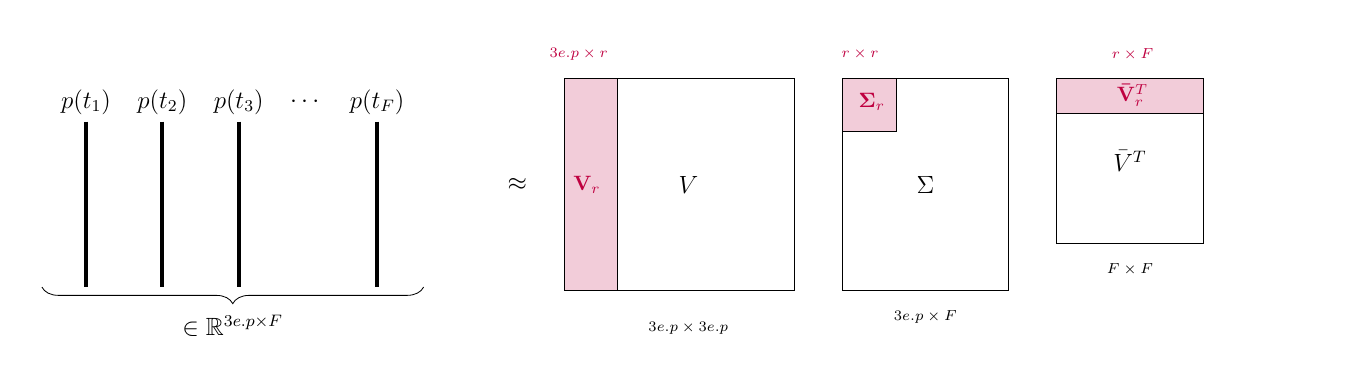}
		\caption{Orthogonal basis computation for vectorized function frames $p(t_i) \in \mathbb{R}^{3N}$ using singular value decomposition. Basis vectors corresponding to the largest singular values are selected.}
		\label{fig:f_vectorized_svd}
	\end{figure*}

	The time evolution of vertex positions $q(t) \in \mathbb{R}^{N \times 3}$, in the absence of damping, is governed by the balance between external forces $f_{\text{ext}}$ and internal forces $f_{\text{in}}(q)$, where the latter is assembled from energy potentials $W_j(q)$ as $f_{\text{in}}(q) = \sum_j W_j(q)$. This yields the standard form of the equation of motion:

	\begin{equation}
		\mathbf{M} \ddot{q} = f_{\text{in}}(q) + f_{\text{ext}},
		\label{eqn:Newtons}
	\end{equation}

	From a variational perspective, advancing the system in time corresponds to solving an energy minimization problem. Specifically, the updated vertex positions $q_{k+1}$ at time step $k+1$ are obtained by solving

	\begin{equation}\label{opt_eq}
		\underset{q_{k+1}}{\text{min}} \ \ \frac{1}{(\delta t)^2} 
		\left\| \mathbf{M}^{\tfrac{1}{2}} ( q_{k+1} - s(q_{k}) ) \right\|^2_F 
		+ \sum_j W_j(q_{k}),
	\end{equation} where $s(q_k)$ denotes the inertial prediction term. This formulation transforms time integration into the solution of an optimization problem, which can be reformulated as a linear system.

	The Projective Dynamics approach \cite{Bouaziz14} leverages this formulation by employing an Alternating Direction Method of Multipliers (ADMM) scheme. It alternates between local constraint projections $p(q_k)$ and a global linear solve that computes the updated positions $q_{k+1}$ given the fixed projections, as summarized in Equation~(\ref{eqn:trajectory}).

	\begin{equation}
		\begin{aligned}
			q_{k+1} &= A^{-1} b_k, \\
			\text{where} \quad & \\
			A &= \frac{1}{(\delta t)^2} \mathbf{M} + \sum_j \omega_j S_j^\top S_j, \\
			b_k &= \frac{1}{(\delta t)^2} \mathbf{M} \ s_k(q_k) + \sum_j S_j^\top p_j(q_k).
		\end{aligned}
		\label{eqn:trajectory}
	\end{equation}

	Let \( C_j \) denote the \( j^{\text{th}} \) constraint manifold. The function \( p_j(q) \) defines the projection of the current positions onto \( C_j \) over a time step \( \delta t \), and \( \omega_j \) is the associated non-negative stiffness weight for the constraint. The projection \( p_j(q) \) is influenced by specific degrees of freedom determined by the selection matrix \( S_j \). Accordingly, the corresponding energy potential takes the form \( W_j = \| S_j q - p_j(q) \| \).

	By defining \( B = \frac{M}{(\delta t)^2} \) and \( u(t) = s_k(q_k) \), the global system in Equation~(\ref{eqn:trajectory}) can be equivalently expressed as
	\[
		A q_{k+1} = B u(t) + \sum_j S_j^\top p_j(q_k),
	\]
	where \( A \) and \( B \) are the system matrices, \( u(t) \) is a time-dependent input term, and \( \sum_j S_j^\top p_j(q_k) \) is a nonlinear correction term updated at each time step for all constraint types, from now on we drop the sum $\sum_j $ and the $j$ subscript for simplicity, and we write the nonlinear term as$S^T p(q(t))$ or $S^T p(t)$ for short.

	Computational complexity for vertex positions can be significantly reduced by constructing a low-dimensional linear subspace \( W \in \mathbb{R}^{N \times k} \), where \( k \ll N \). The full solution \( q \) is then approximated as \( q \approx W \hat{q} \), with \( \hat{q} \in \mathbb{R}^{k \times 3} \) representing the reduced coordinates. This allows the system to be solved in the low-dimensional subspace, avoiding the high computational cost of solving the full system in Equation~(\ref{eqn:trajectory}). The resulting reduced system takes the form:

	\begin{equation}
		\label{eqn:reducedSys}
		\tilde{A} \hat{q}_{k+1} = \tilde{b}_k,
	\end{equation} where the reduced matrices \( \tilde{A} \) and \( \tilde{b}_k \) are derived via projection.

	The basis \( W \) can be constructed using linear blend skinning subspaces \cite{Brandt18}, or from snapshot-based methods that capture representative solutions from the full space \cite{Monem23}.

	However, when the nonlinear function \( p : \mathbb{R}^{N \times 3} \rightarrow \mathbb{R}^{e \cdot \tilde{p} \times 3} \) depends on constrained elements such as vertices, edges, faces, or tetrahedra, the computational cost of evaluating the nonlinear term in Equation~(\ref{eqn:reducedSys}) remains high. In many cases, the number of constrained elements \( e \) is larger than the number of vertices \( N \), and \( \tilde{p} \) represents the row dimension of the matrix-valued function \( p(\cdot) \). As a result, the evaluation of \( p(W \hat{q}(t)) \) can still dominate the overall computation, despite the reduced dimensionality of the vertex space.

	\begin{equation}
		\underbrace{W^T \ E \ W}_{k \times k \times 3} \hat{q}(t) =  \underbrace{W^T \ B W}_{k \times k \times 3} \hat{u}(t) + \underbrace{W^T}_{k \times N \times 3} \overbrace{S^T}^{N \times e.p} \ \underbrace{p(W \ \hat{q}(t))}_{e.p \times 3}.
		\label{eqn:reducedSys}
	\end{equation}
	
	Basis computed using snapshot-based methods \cite{Monem23, Monem25} have demonstrated improved rotational expressiveness and numerical stability compared to linear blend skinning subspaces \cite{Brandt18} when applied to vertex position reduction. Motivated by the success of snapshot-based subspaces, the Discrete Empirical Interpolation Method (DEIM) emerges as a promising candidate for reducing the computational complexity of evaluating the nonlinear term \( p(W \hat{q}(t)) \). This study investigates the compatibility of DEIM for reducing constraint projections within the projective dynamics framework.

	The DEIM method consists of two main steps \cite{Chaturantabut10}. First, \( F \) flattened time instances (snapshots) of the nonlinear function \( f \in \mathbb{R}^{3p} \) are collected. Singular value decomposition (SVD) is then applied to the resulting snapshot matrix. The linear vector basis of rank \( r \) that best represents the snapshots corresponds to the first \( r \) columns of the basis matrix, denoted \( V = V_r \), as illustrated in Figure~\ref{fig:f_vectorized_svd}. The basis matrix \( \left[ v_1, v_2, \dots, v_r \right] \) is then used for the reduction process. 

	For simplicity, we omit the subscript and use \( V \) to denote the basis matrix throughout. The nonlinear function is then approximated as
	\[
		p(t) \approx V \hat{p}(t),
	\]
	where \( \hat{p}(t) \in \mathbb{R}^{r \times 3} \) is the latent representation of the nonlinear function.

	In each iteration of DEIM, the selection matrix \( P^\top \) is updated by adding a standard basis vector \( \wp_l \in \mathbb{R}^{3ep} \), which selects the index corresponding to the row of the flattened vector basis \( v_i \in \mathbb{R}^{3ep} \) that exhibits the largest projection error relative to the subspace spanned by the previously selected basis vectors \( \left[ v_l \right]_{1}^{i-1} \). The core idea is to reduce the overall approximation error by enforcing exactness at the points with the largest residuals.

\begin{figure*}[h!]
		\centering
		\includegraphics[width=\textwidth, height=0.5\textwidth, keepaspectratio]{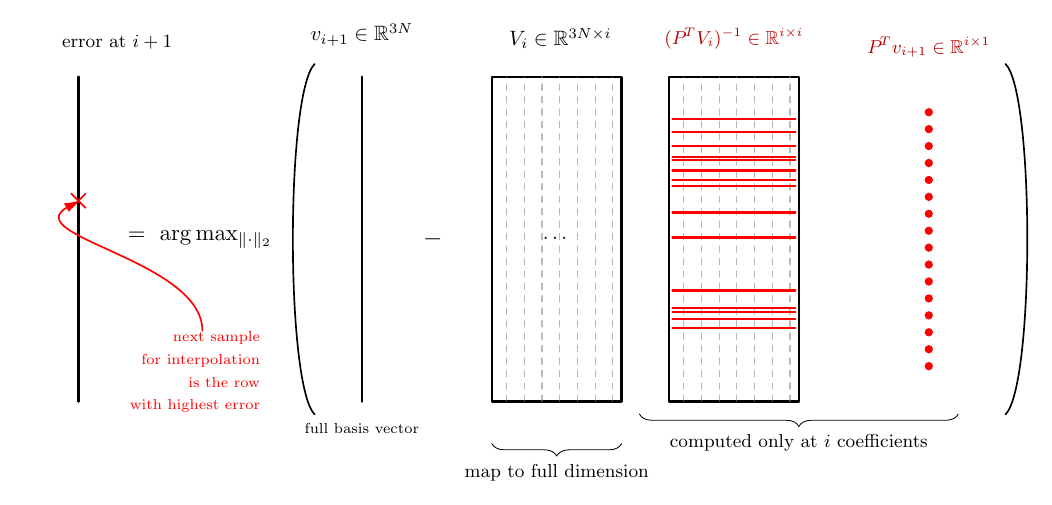}
		\caption{For each vector basis  $ v_{i+1},  \ 1 \leq i \leq r$ The point of the highest interpolation error is determined through DEIM solve and added to the selection matrix $P^T$ for the next iteration.}
		\label{fig:deim_solve_visualized}
	\end{figure*}

	In the second step of DEIM, the \( r \) points that are least represented by the current vector basis are gradually identified, and interpolation equality is enforced at these points. Specifically, for each basis vector \( v_i \), the condition \( P^\top V v_i = P^\top v_i \) is enforced, where the selection matrix \( P^\top \) defines the chosen interpolation points. This process is illustrated in Figure~\ref{fig:deim_solve_visualized} and described in detail in Algorithm \ref{alg:DEIM}.

	This step enhances the accuracy of basis interpolation and provides a principled way to select a small set of mesh points that can be used for online evaluation of nonlinear forces. In other words, the interpolation points selected by DEIM are used during real-time computation to approximate the nonlinear function across the entire mesh, as described in Equation~(\ref{eqn:deim_approximate}).

	The solver \( (P^\top V)^{-1} \in \mathbb{R}^{r \times r \times 3} \) is computed once offline. During the online (real-time) phase, only the values \( P^\top p(t) \in \mathbb{R}^{r \times 3} \) need to be evaluated, significantly reducing computational cost.

	\begin{equation}
		V \ \hat{p}(t) = V \ (P^\top V)^{-1} \ P^\top \ p(t),
		\label{eqn:deim_approximate}
	\end{equation}

	% visualize original deim solve

	% DEIM original algorithm
 	\begin{algorithm}[]
    	\textbf{Input:} $\left[ v_l \right]_{l=1}^r \subset \mathbb{R}^{3e.p}$ linear independent basis vectors associated to vectorized simulations.
        \newline
        
        $\left[ \| \rho\|, \wp_1 \right] = \max(\| V_1\|)$\\
        $ V = \left[ v_1 \right], \quad P = \left[ \wp_1 \right]$  \\
        \For{$ 2 \leq l \leq m$ :}{
            $(P^T V) c = P^T v_l$ for $c$\\
            $r = v_l - U c$\\
            $\left[ \| \rho\|, \wp_l \right] = \max(\| r\|)$\\
            $V \leftarrow \left[ V, v_l\right] , \quad P \leftarrow  \left[ P, \wp_l \right]$
        } 
        \caption{DEIM}
    	\label{alg:DEIM}
    \end{algorithm}

	The core idea of DEIM is to evaluate the pointwise nonlinear function \( p(\cdot) \) only at a small number of rows selected by the matrix \( P^\top \), and to use a compact solver \( \left( P^\top V \right)^{-1} \in \mathbb{R}^{r \times r} \) to compute the latent representation of the nonlinear force. Traditionally, DEIM has been applied to vector-valued nonlinear functions that both operate on and return values in the same state-space manifold—for example, vertex positions.

	In the case of vertex-based constraints, such as bending, we have \( \tilde{p} = 1 \) and \( e = N \), and the nonlinear projection becomes a vector-valued function \( p : \mathbb{R}^{N \times 3} \rightarrow \mathbb{R}^{N \times 3} \).

	However, when attempting to apply DEIM to a broader range of constraint projections in physics-based animation, a key challenge arises: the nonlinear projections \( p_j(q(t)) \) typically take inputs from the position space, but do not necessarily yield outputs in the same space. More generally, for a given number of constrained elements \( e \), the projection onto the constraint manifold \( C_j \) is represented by a matrix-valued function \( p_j(q(t)) \in \mathbb{R}^{e \cdot \tilde{p}_j \times 3} \), where \( \tilde{p}_j \) denotes the row dimension of the projection per element.

	For instance, a bending constraint is computed on vertices, resulting in a \( 1 \times 3 \) vector per element, with \( \tilde{p} = 1 \) and \( e_{\text{vertices}} = N \). In contrast, strain forces may be computed as \( 2 \times 3 \) vectors over \( e_{\text{faces}} \) elements, with \( \tilde{p} = 2 \), and so on for other constraint types.

	In DEIM, the linearly independent vector basis is computed from flattened snapshots. Identifying interpolation points in this way results in separate indices for each spatial dimension (\( x \), \( y \), \( z \)), which is not necessarily suitable for constraint projection, as it disrupts the inherent three-dimensional correlation between mesh elements.

\begin{figure}[h!]
		\centering
		\begin{adjustbox}{width=0.8\textwidth} % or =\textwidth, =0.8\textwidth, etc.
		\begin{minipage}{\textwidth}

		\includegraphics[width=\textwidth, height=0.5\textwidth, keepaspectratio]{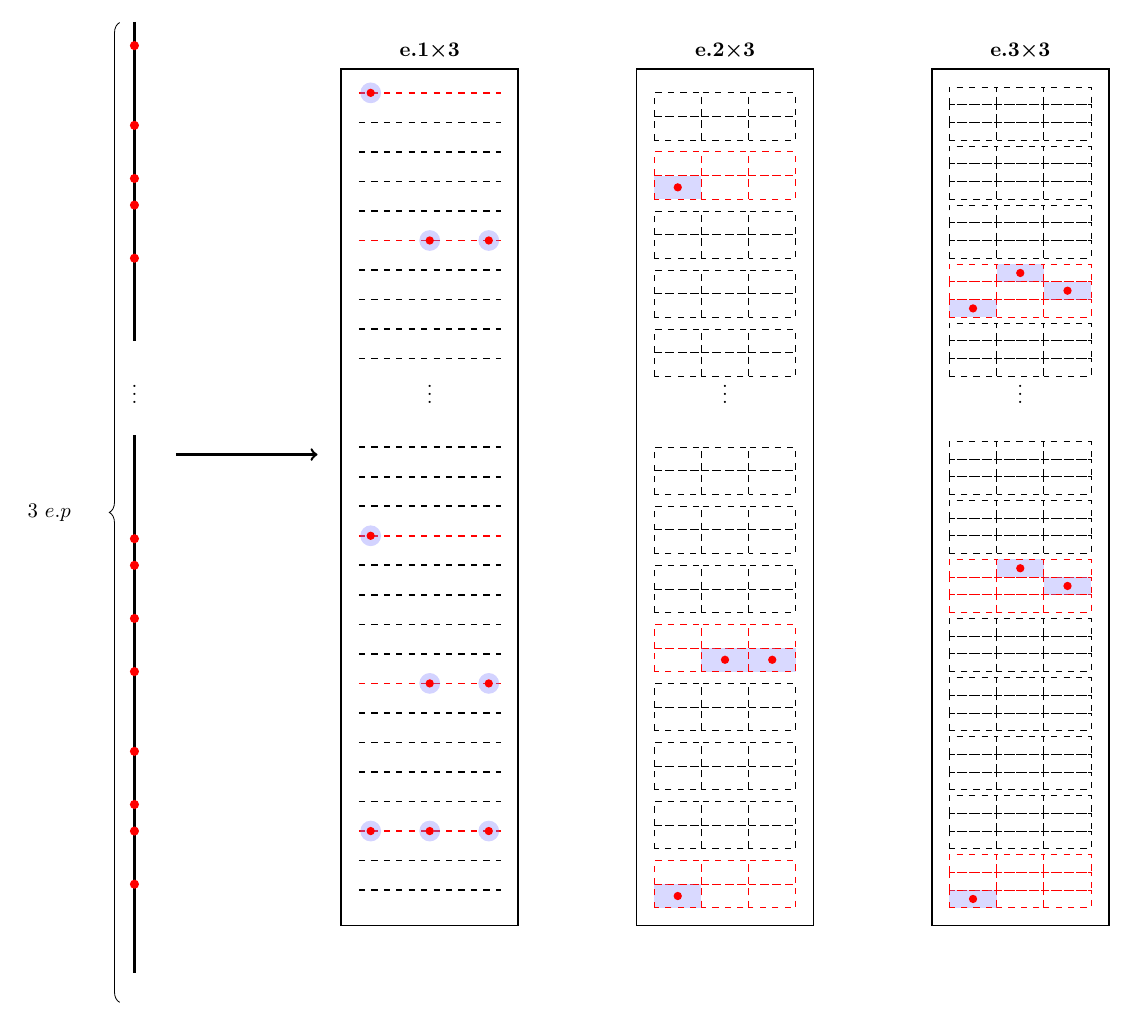}
		\end{minipage}
		\end{adjustbox}

		\caption{Vectorizing constraint projections interpolation points determinations is inconvenient and inefficient.}
		\label{fig:deim_vectorized_sampling}
	\end{figure}

	If projection snapshots are flattened into vectors in \( \mathbb{R}^{3e \cdot \tilde{p}_j} \) for basis computation, and DEIM is applied to select the \( r \) rows with the largest error, the resulting indices often correspond to scattered entries from different projections and different elements. This behavior is visually illustrated in Figure~\ref{fig:deim_vectorized_sampling}. In practice, this means that for each selected entry, one must compute the full \( \tilde{p} \times 3 \) projection in order to use only a single component—leading to unnecessary computation that contradicts the purpose of reduced models. 

	Alternatively, one could implement nested conditional logic to extract only the necessary entries for each case, but this significantly increases code complexity and reduces generality. For instance, applying DEIM to tetrahedral strain constraints would require computing the full \( 3 \times 3 \) deformation matrix for many elements, only to extract a single value from the nine possible entries. This not only wastes resources but also results in poor runtime efficiency. 

	In deformable character physics simulations, it is both intuitive and essential to preserve the full 3D correlation during interpolation. Maintaining this coherence is critical to ensuring accurate and efficient approximation of nonlinear constraint projections.

	However, in our opinion, DEIM methods in their current form are not well-suited for selecting a small set of constrained elements for force interpolation. This limitation arises primarily from the fact that constrained elements are not always vertices in the state space; they can also include edges, triangles, or tetrahedra. Moreover, standard DEIM does not preserve the three-dimensional correlation between spatial components, which is essential in computer graphics applications.

	% DEIM vectorized sampling.
	
	In the following section, we present experimental examples to evaluate the capability of DEIM in selecting a small set of constrained elements for constructing a latent subspace suitable for nonlinear constraint projections.

%%%%%%%%%%%%%%%%%%%%%%%%%%%%%%%%%%%%%%%%%%%%%%%%%%%%%%%%%%%%%%%%%%%%%%%%%%%%%%%%
%%%%%%%%%%%%%%%%%%%%%%%%%%%%%%%%%%%%%%%%%%%%%%%%%%%%%%%%%%%%%%%%%%%%%%%%%%%%%%%%
\section{Implementations and Discussion}

	To evaluate the capability of the DEIM algorithm in selecting a small subset of elements for computing reduced nonlinear forces, we conducted a study on two simple examples of deformable meshes, each subjected to different types of constraint projections. Specifically, we tested the method on three constraint types: bending, spring, and strain.

	For each example, \( F \) frames were collected, where each frame corresponds to a constraint force snapshot \( f \in \mathbb{R}^{e \cdot \tilde{p} \times 3} \), with \( 1 \leq f \leq F \). These snapshots were reshaped into a matrix of size \( \mathbb{R}^{3e \cdot \tilde{p} \times F} \) to compute a POD basis using singular value decomposition (SVD). This process yields at most \( r = F \) orthonormal basis vectors, and the DEIM selection matrix \( P^\top \) can select up to \( F \) row indices.

	As a result, the computed basis has the form \( V \in \mathbb{R}^{3e \cdot \tilde{p} \times r} \), which can be separated into dimension-specific components \( V_x, V_y, V_z \in \mathbb{R}^{e \cdot \tilde{p} \times r} \) for each of the \( x \), \( y \), and \( z \) spatial dimensions.

	To determine the DEIM interpolation points, we run Algorithm~\ref{alg:DEIM} with a key modification: the three-dimensional correlation is preserved throughout the process. As a result, the basis vectors are treated as blocks in \( \left[ v_l \right]_{l=1}^r \subset \mathbb{R}^{e \cdot \tilde{p} \times 3} \).

	At each iteration, the linear system \( (P^\top V) c = P^\top v_l \) is solved for \( c \in \mathbb{R}^{r \times 3} \), and the interpolation error \( r = v_l - V c \) is computed jointly across all three spatial dimensions (\( x, y, z \)). This ensures that the interpolation error is evaluated with the full 3D structure of the basis preserved.

	In Step 7 of the algorithm, we identify the row index in the range \( 1 \leq i \leq e \cdot \tilde{p} \) with the largest \( \| \cdot \|_2 \) norm of the residual and update the selection matrix \( P^\top \) accordingly.
	The goal of this experiment is to evaluate the compatibility of the proposed DEIM-based technique for computing suitable reduced bases in the context of physics-based animation, rather than to challenge its ability to represent complex dynamics. To this end, we designed two basic experiments.

	In the first test, we simulated non-reduced cloth dynamics under three types of constraint projections: bending, spring, and strain. For each constraint manifold \( C_j \), we collected snapshots \( p_j(q(t)) \), using 25, 100, and 100 frames respectively. Reduced bases were then computed as described in the previous section.

\begin{figure}[h!]
	\centering
	% tweak horizontal gaps between panels (column padding)
	\setlength{\tabcolsep}{3pt}
	% remove extra row height in the tabular
	\renewcommand{\arraystretch}{0}

	% inside a figure or figure* environment
	\centering
	\begin{adjustbox}{width=0.5\textwidth} % or =\textwidth, =0.8\textwidth, etc.
		\begin{minipage}{0.80\textwidth}
	% ---------- Row 1 ----------
	\begin{subfigure}[t]{\textwidth}
	\centering
	\includegraphics[trim={21cm 6cm 21cm 4cm},clip,width=0.31\textwidth]{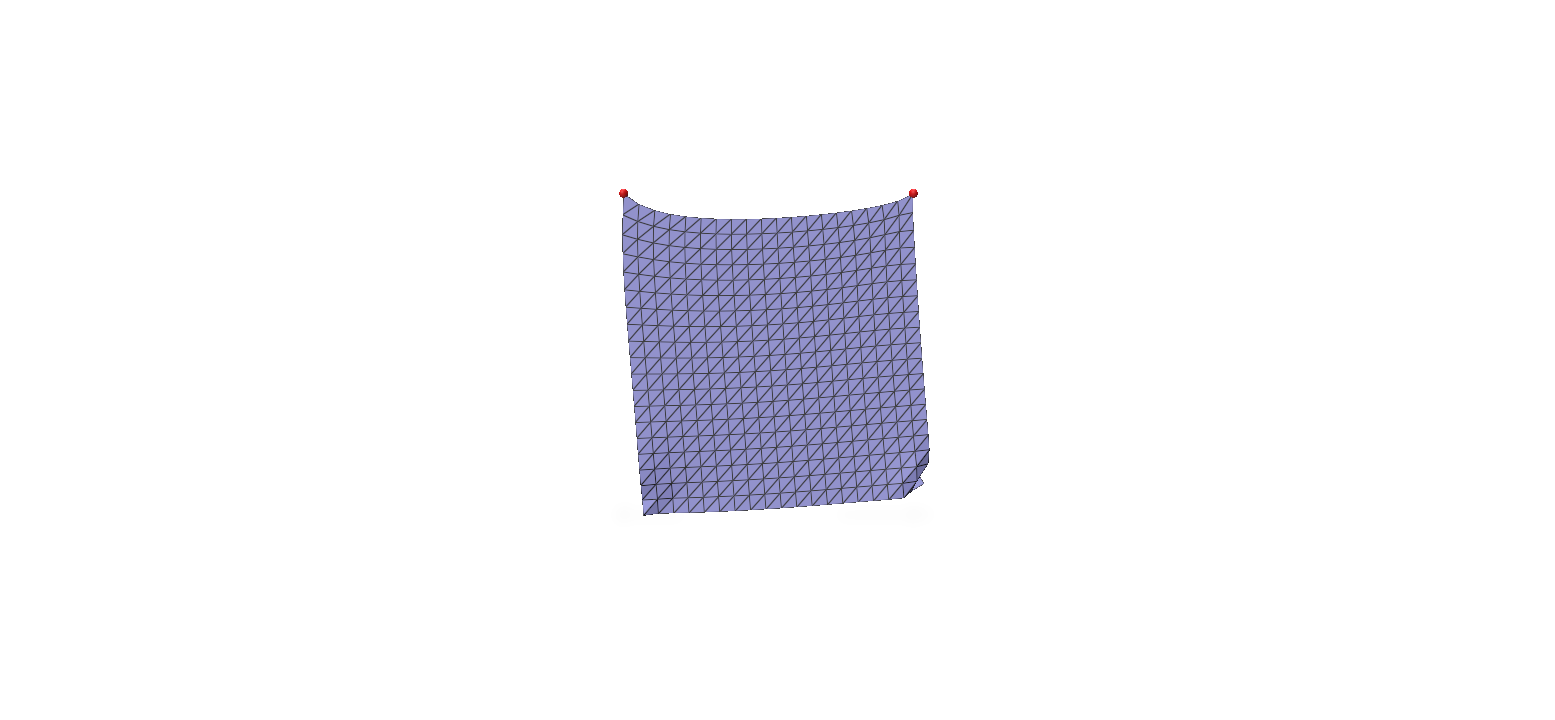}\hfill
	\includegraphics[trim={21cm 6cm 21cm 4cm},clip,width=0.31\textwidth]{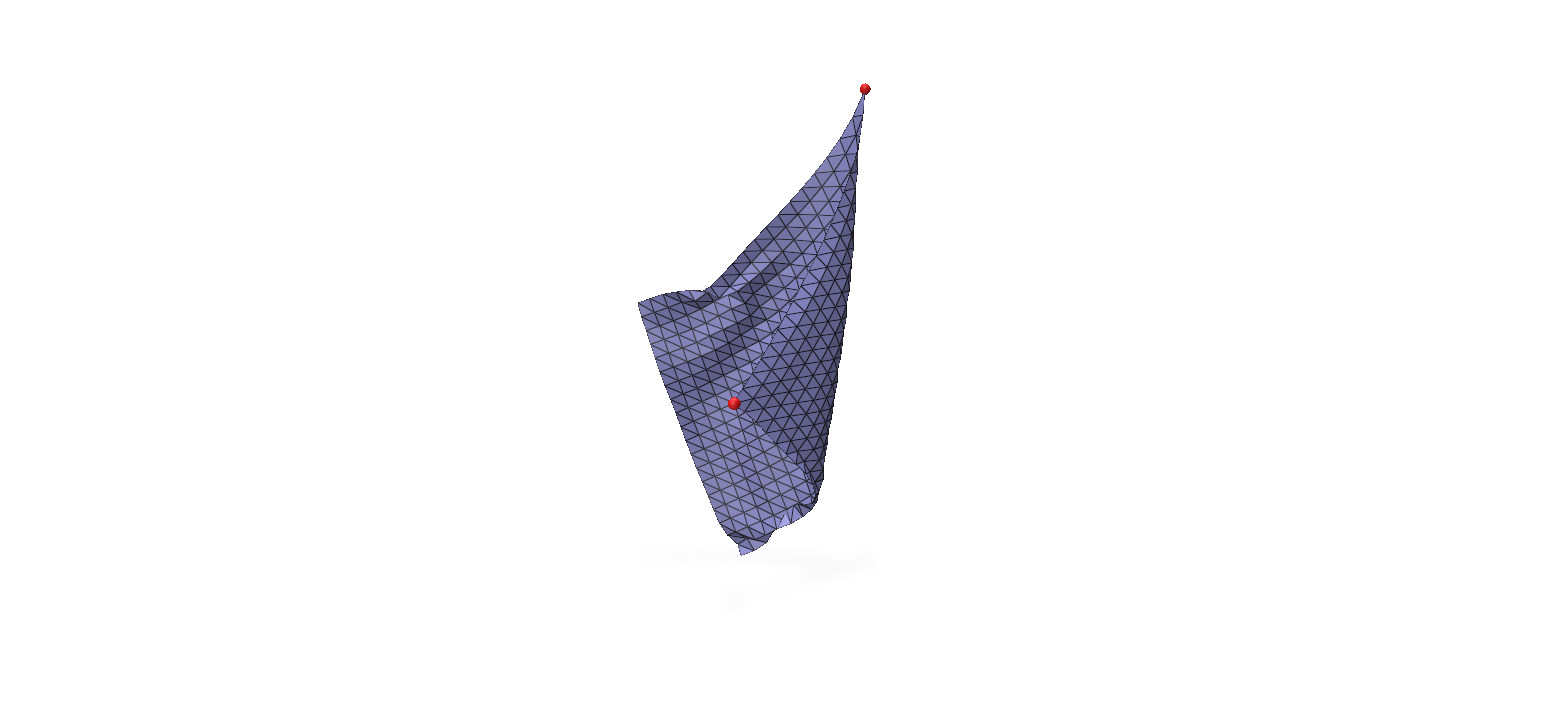}\hfill
	\includegraphics[trim={21cm 6cm 21cm 4cm},clip,width=0.31\textwidth]{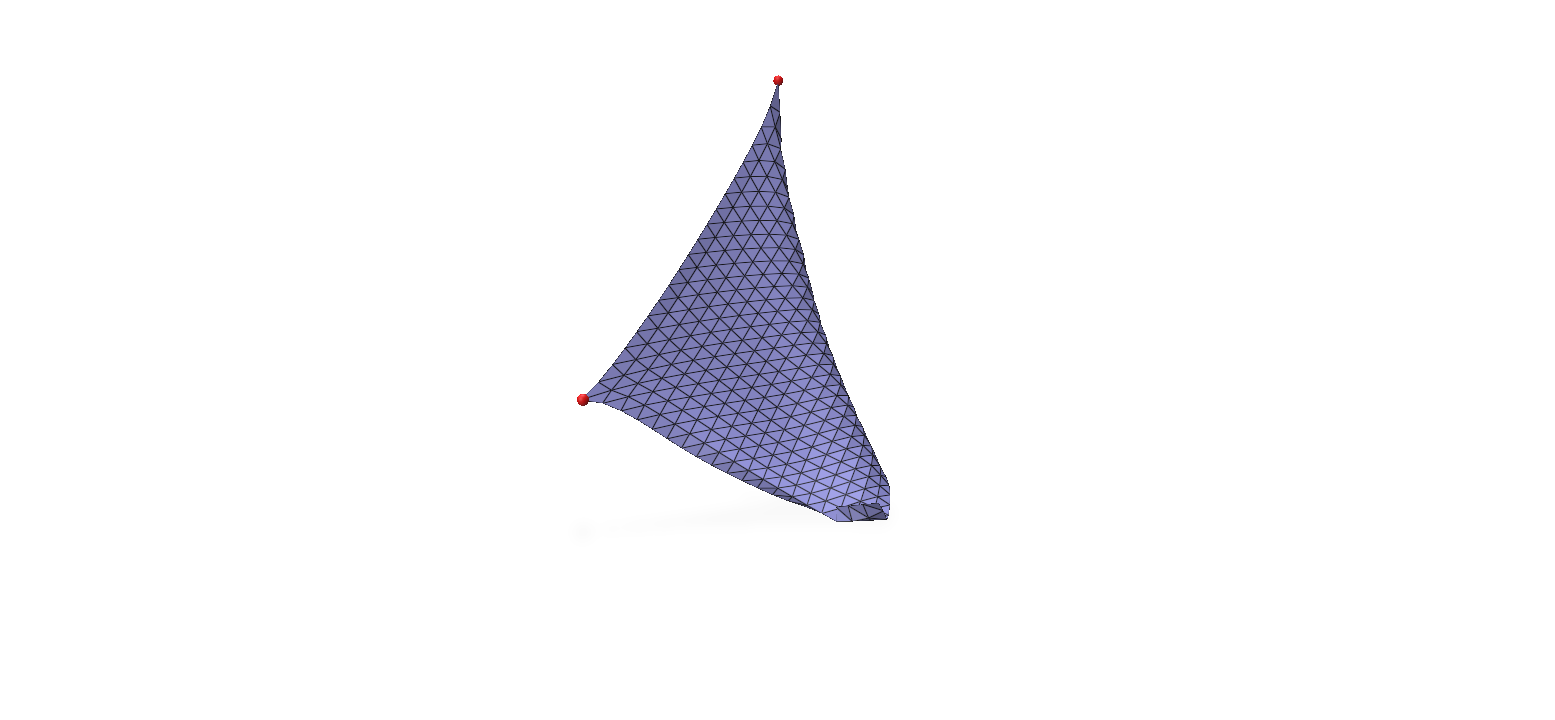}
	\caption{Non reduced simulations.}
	\end{subfigure}

	\vspace{4pt}

	% ---------- Row 2 ----------
	\begin{subfigure}[t]{\textwidth}
	\centering
	\includegraphics[trim={19cm 2cm 19cm 4cm},clip,width=0.31\textwidth]{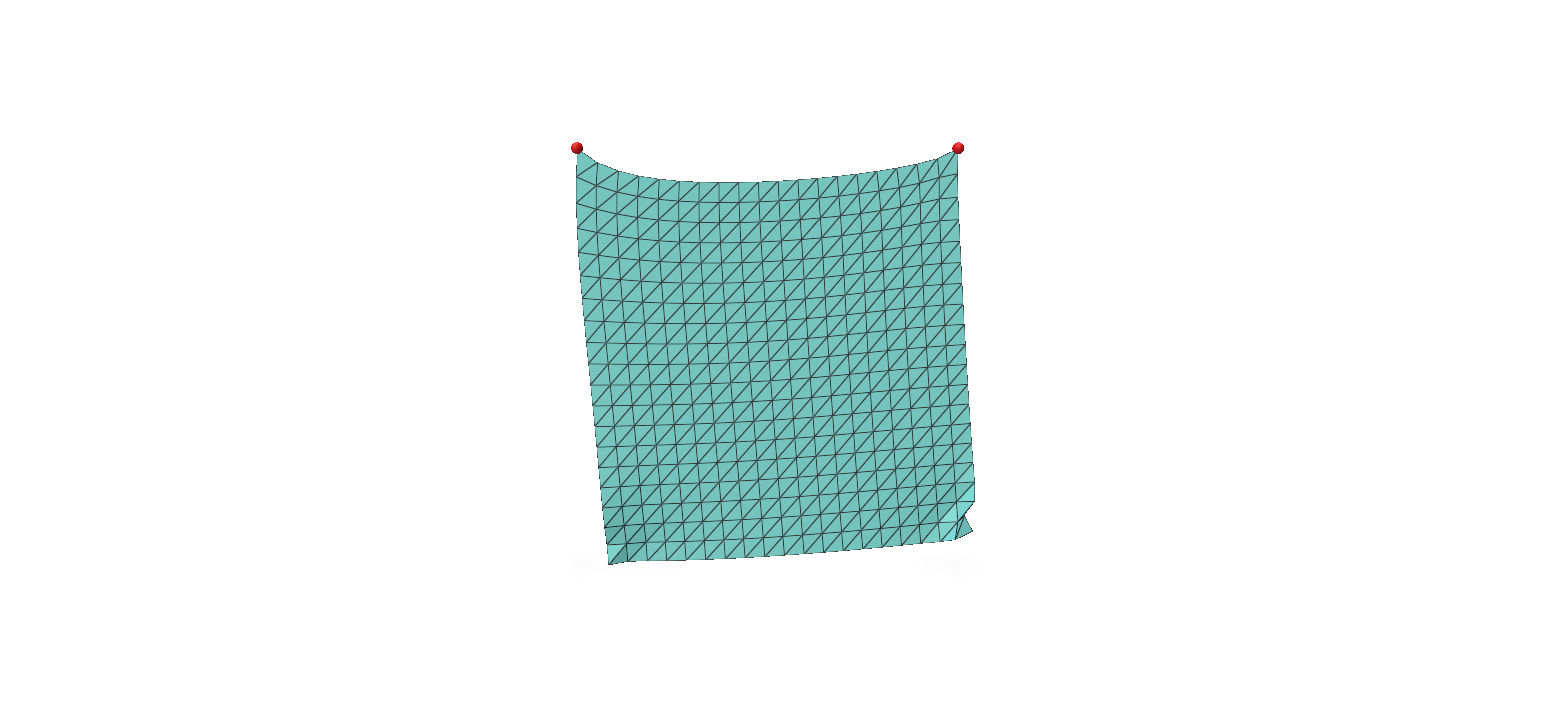}\hfill
	\includegraphics[trim={19cm 2cm 19cm 4cm},clip,width=0.31\textwidth]{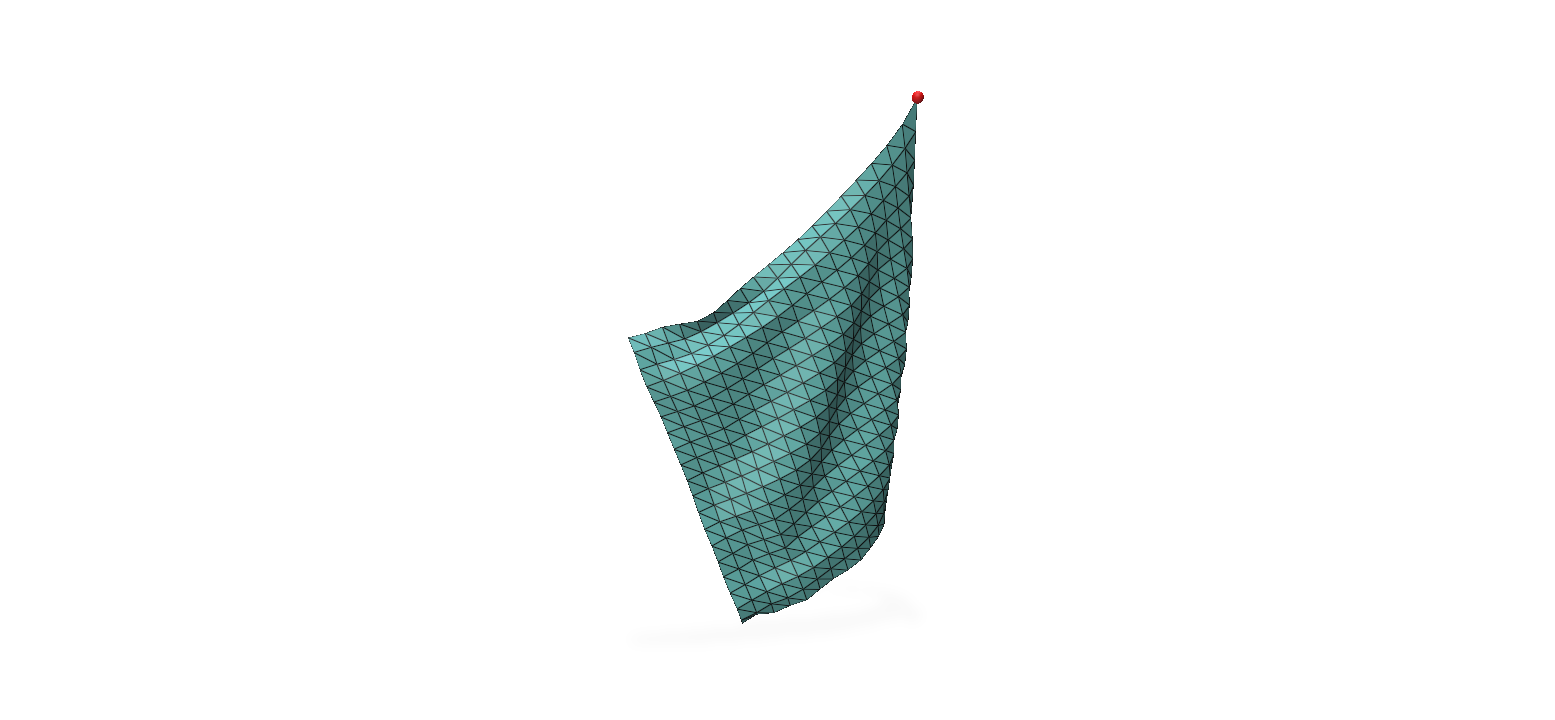}\hfill
	\includegraphics[trim={20cm 1cm 20cm 4cm},clip,width=0.31\textwidth]{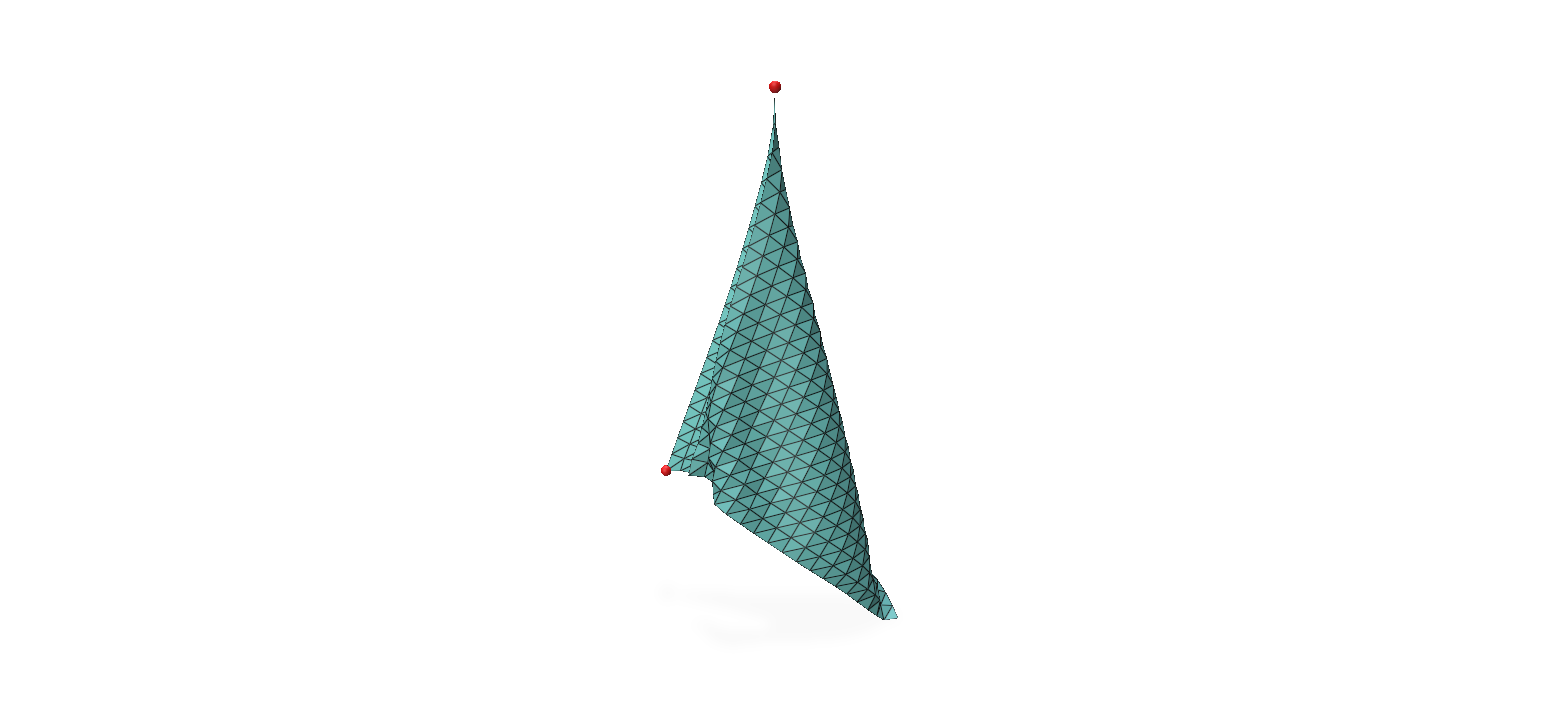}
	\caption{Reduced vertex bending only.}
	\end{subfigure}

	\vspace{4pt}

	% ---------- Row 3 ----------
	\begin{subfigure}[t]{\textwidth}
	\centering
	\includegraphics[trim={20cm 4cm 20cm 2cm},clip,width=0.31\textwidth]{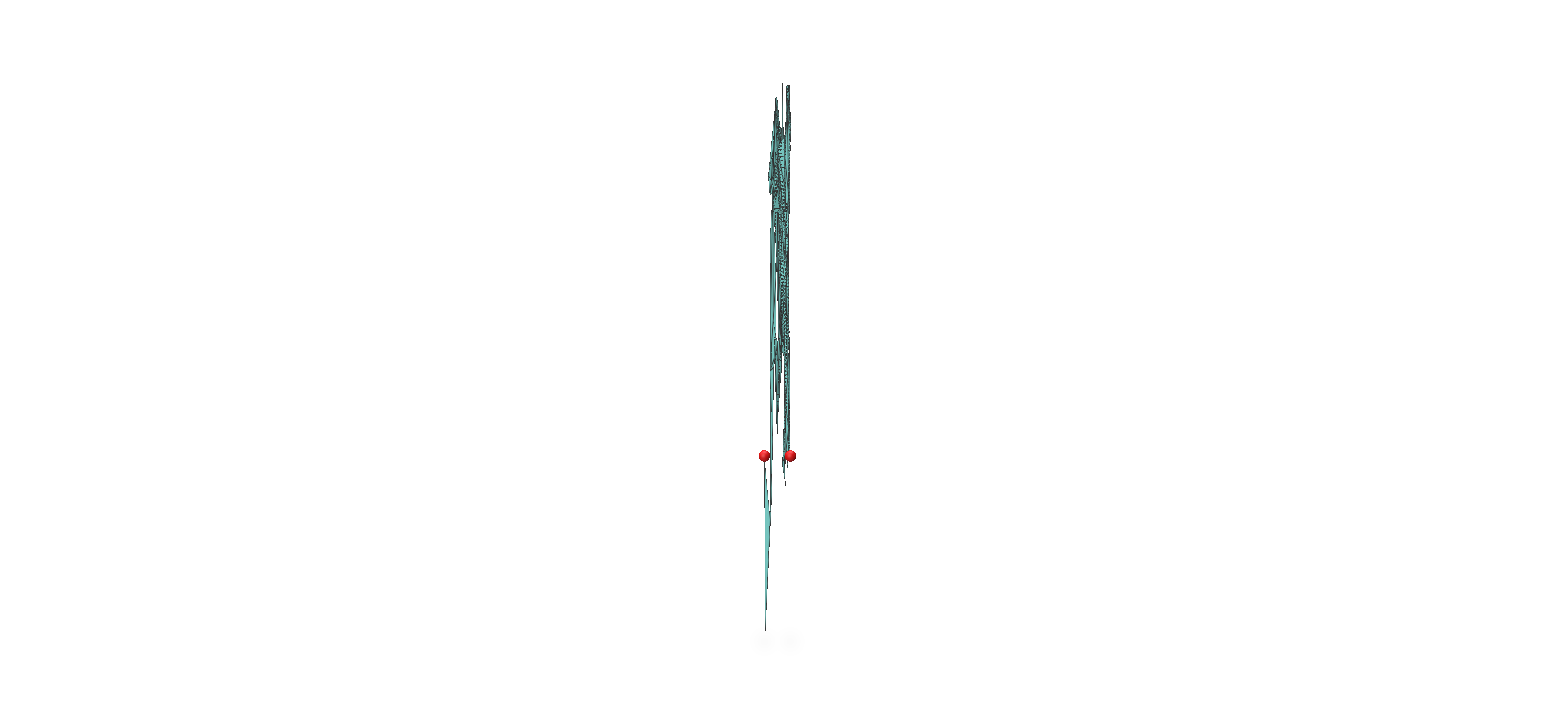}\hfill
	\includegraphics[trim={20cm 4cm 20cm 2cm},clip,width=0.31\textwidth]{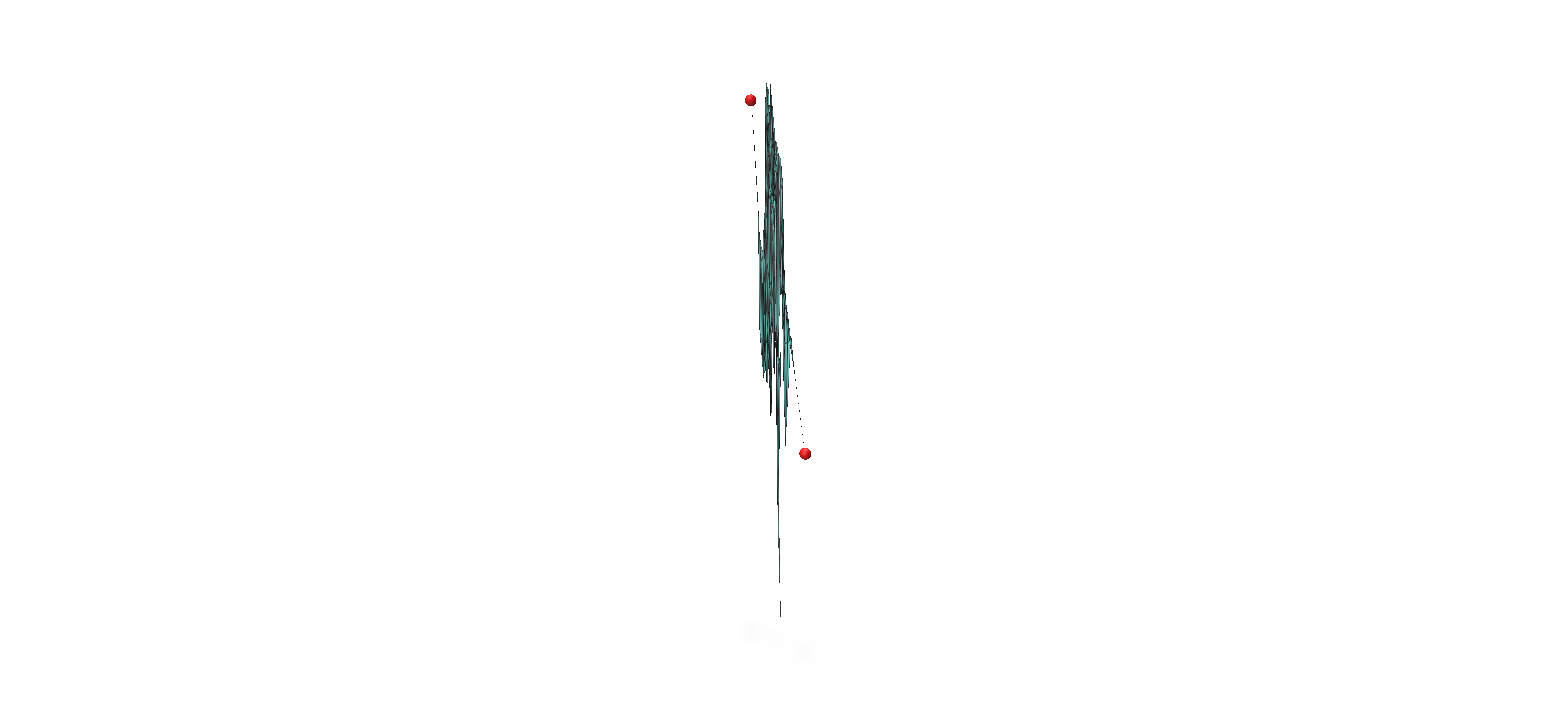}\hfill
	\includegraphics[trim={20cm 4cm 20cm 2cm},clip,width=0.31\textwidth]{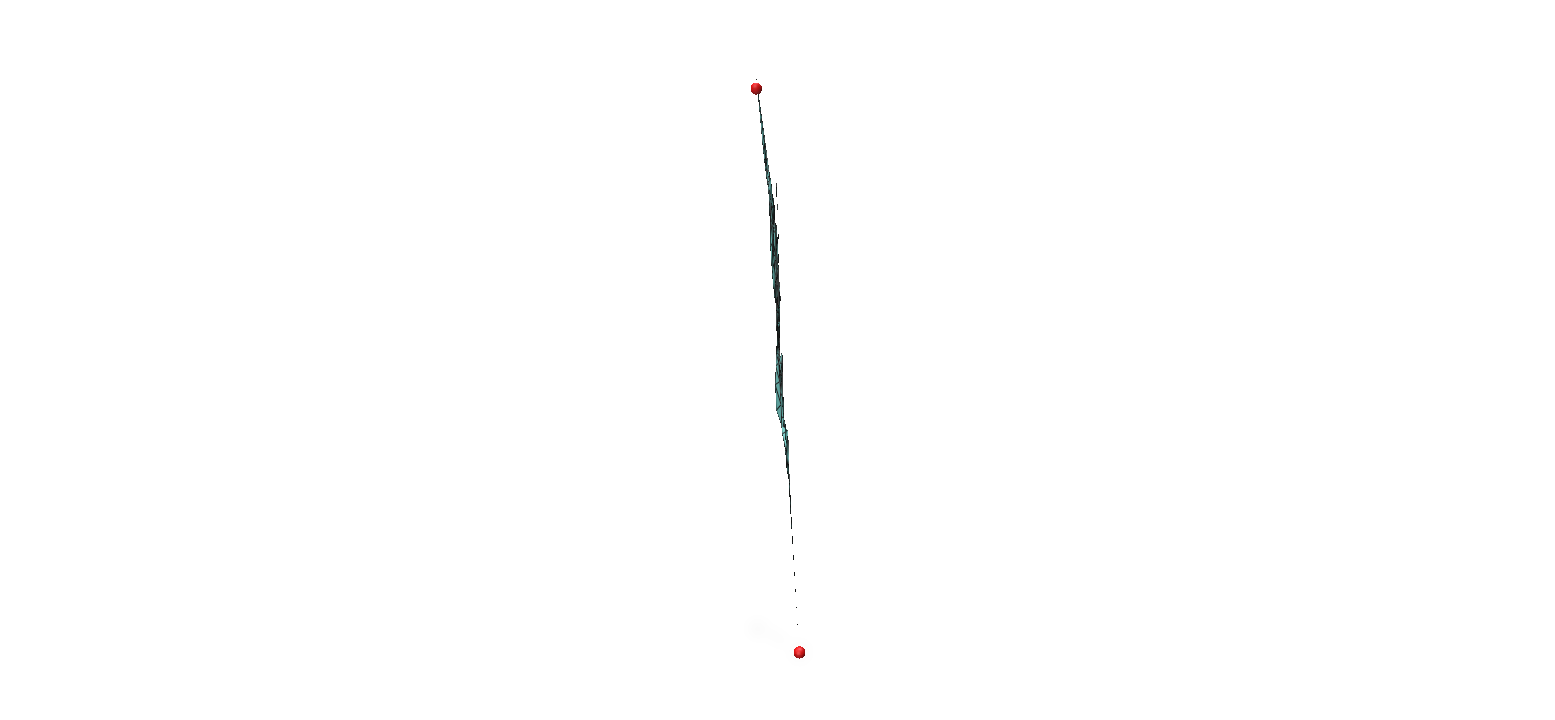}
	\caption{Reduced edge spring only.}
	\end{subfigure}

	% ---------- Row 4 ----------
	\begin{subfigure}[t]{\textwidth}
	\centering
	\includegraphics[trim={19cm 4cm 19cm 2cm},clip,width=0.31\textwidth]{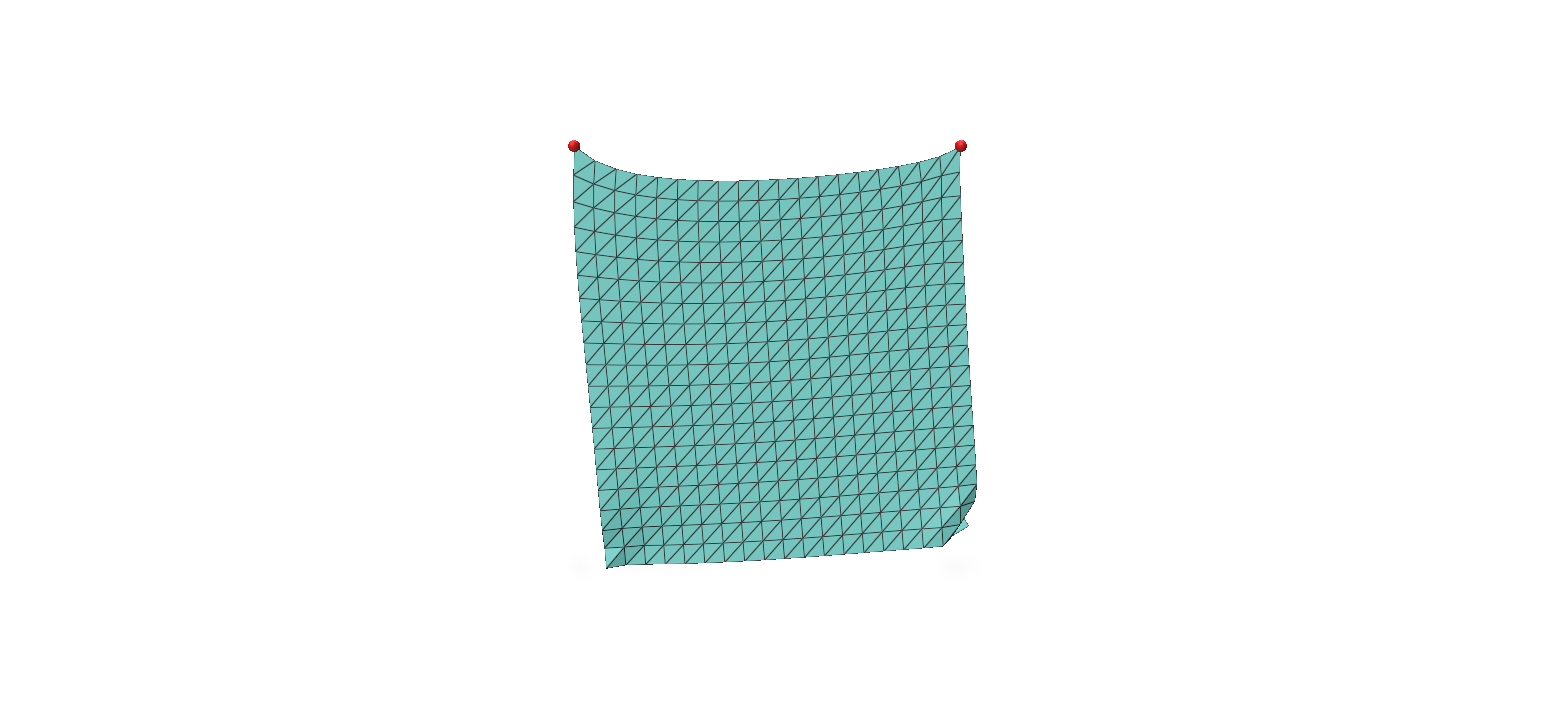}\hfill
	\includegraphics[trim={19cm 4cm 19cm 2cm},clip,width=0.31\textwidth]{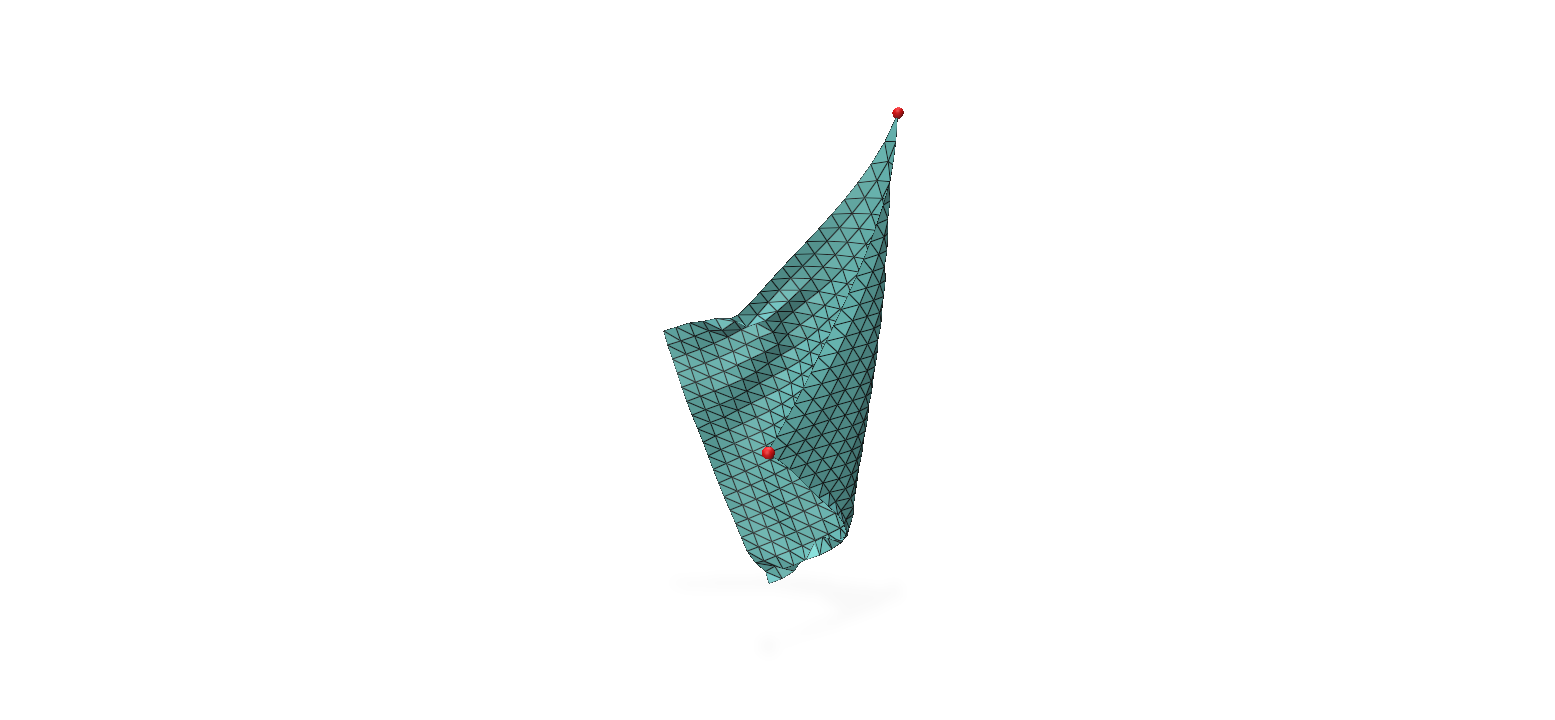}\hfill
	\includegraphics[trim={19cm 4cm 19cm 2cm},clip,width=0.31\textwidth]{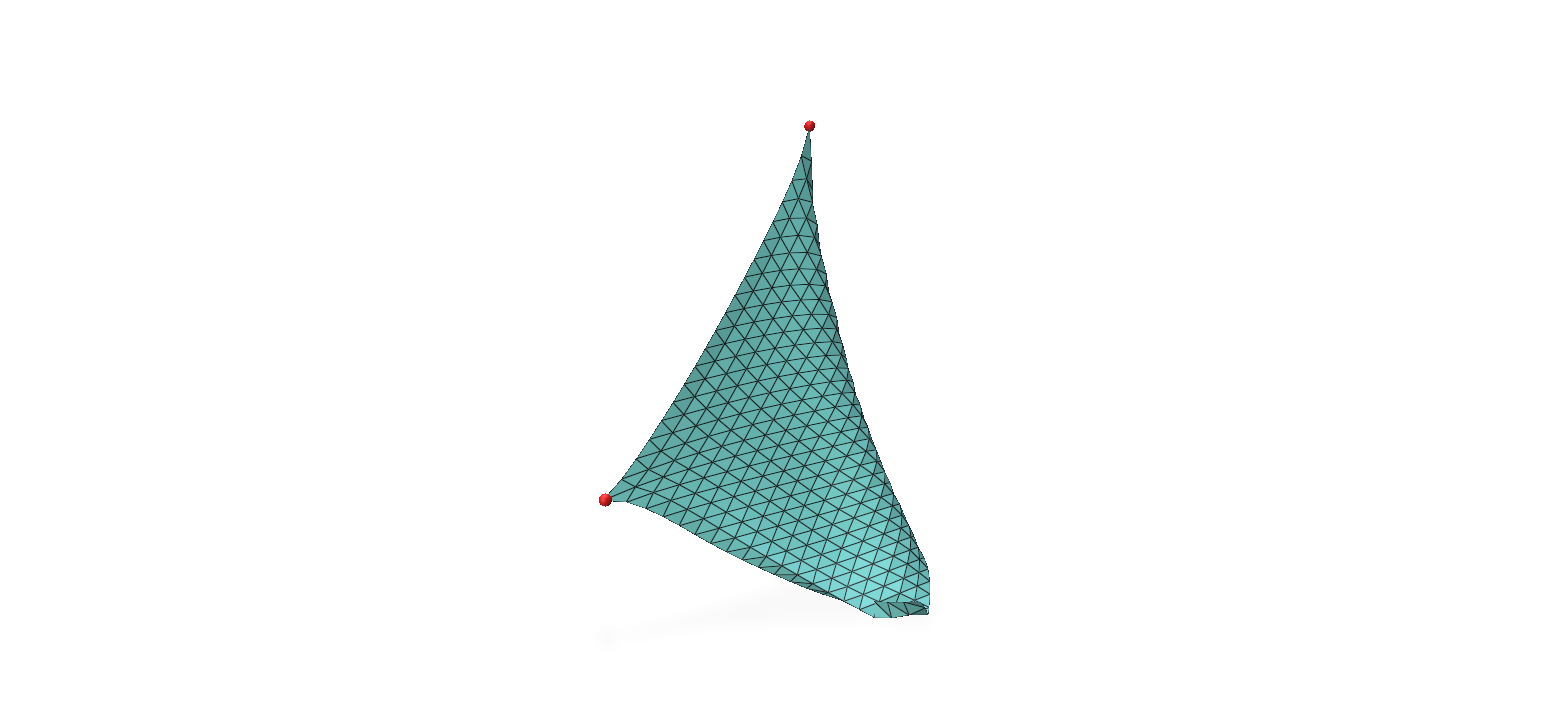}
	\caption{Reduced edge spring only.}
	\end{subfigure}

		\end{minipage}
	\end{adjustbox}

	\caption{Misleading experiment where basis where computed on fully constrained mesh with three types of constraints together. Reconstruction relative error measured while increasing the basis carinality, on the cloth simulation training frames for different constraint projections..Cloth mesh simulation shown, under three sets of constraints: vertex bending with homogenous weights $w= 10^{-1}$, edges spring with $w= 10^{6}$, and triangle strain at $w= 10^{-2}$. Comparing between (a) full non reduced constraint projections (b) only reduced vertex bending projections and (c) only reduced edge spring projections.}
	\label{fig:cloth_sim_grid}
	\end{figure}

	To isolate the effect of each constraint's reduction, we visualized simulations where only one constraint type was reduced at a time, while the others were kept in full dimension. We observed that when only edge spring forces were reduced and the rest ran at full resolution, the simulation became completely distorted. In contrast, reducing only the bending or only the strain projections resulted in visually plausible simulations, albeit with some minor artifacts (see Figure~\ref{fig:cloth_sim_grid}).

\begin{figure*}[h!]

		\centering
		% tweak horizontal gaps between panels (column padding)
		\setlength{\tabcolsep}{3pt}
		% remove extra row height in the tabular
		\renewcommand{\arraystretch}{0}

		% inside a figure or figure* environment
		\centering
		\begin{adjustbox}{width=\textwidth} % or =\textwidth, =0.8\textwidth, etc.
		\begin{minipage}{\textwidth}

			% ---------- Row 1 ----------
		\begin{subfigure}[t]{\textwidth}
		\centering
		\includegraphics[trim={12cm 8cm 12cm 8cm},clip,width=0.4\textwidth]{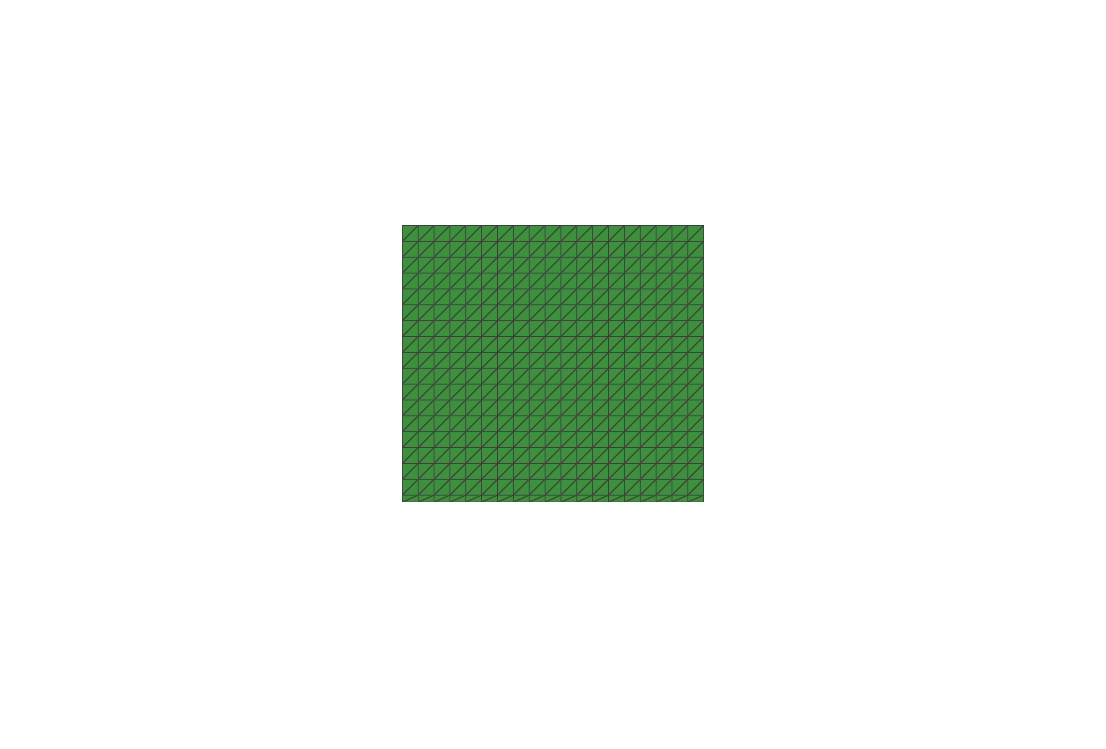}\hfill
		\includegraphics[trim={12cm 8cm 12cm 8cm},clip,width=0.4\textwidth]{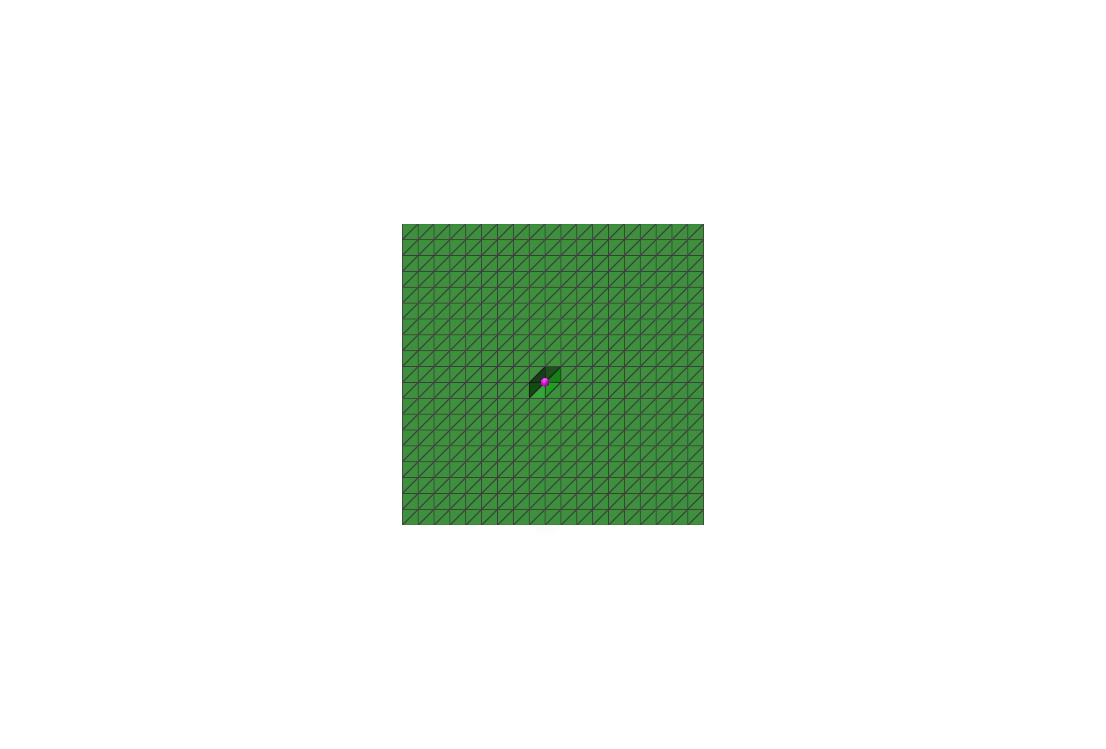}
		
		\caption{FOM}
		\end{subfigure}
		\vspace{-4pt}
		% ---------- Row 2 ----------
		\begin{subfigure}[t]{\textwidth}
		\centering
		\newcommand{\wNine}{\dimexpr\linewidth/9\relax}
\begin{minipage}{0.45\textwidth}
		\blankgraphics{2\wNine}{2.8\wNine}%
		\includegraphics[trim={12cm 5cm 12cm 5cm},clip,width=0.32\textwidth]{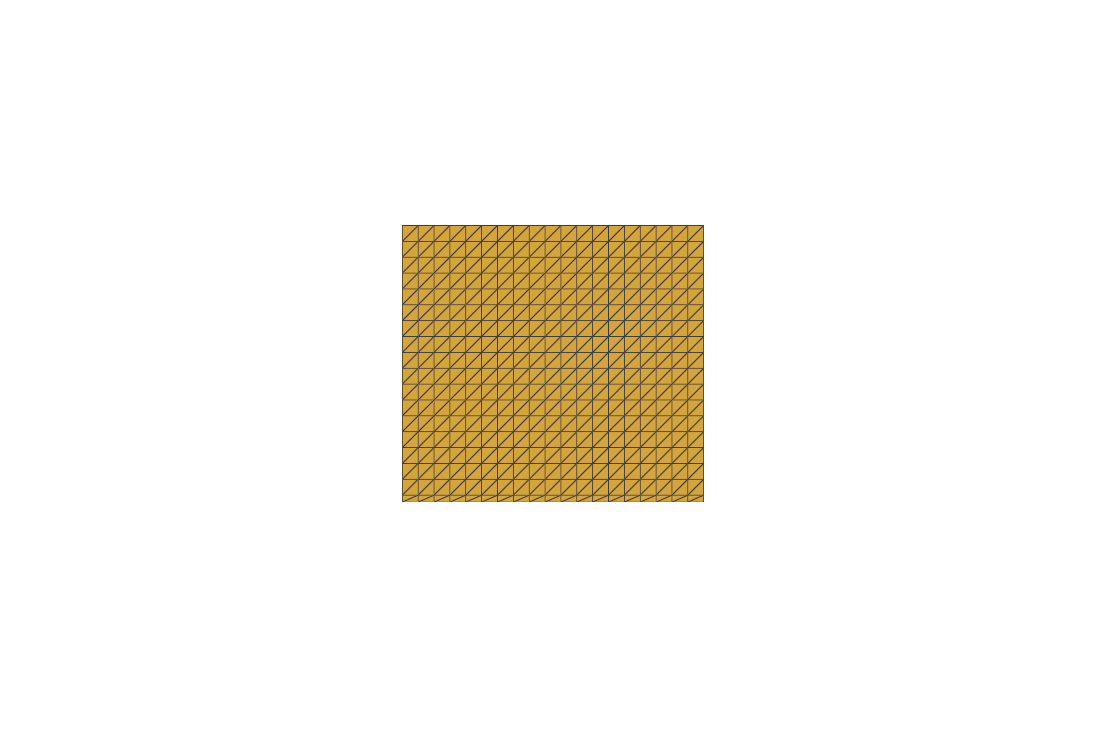}%
		\blankgraphics{2\wNine}{2.8\wNine}%
\end{minipage}
		\hfill
\begin{minipage}{0.45\textwidth}
		\includegraphics[trim={12cm 5cm 12cm 5cm},clip,width=0.32\textwidth]{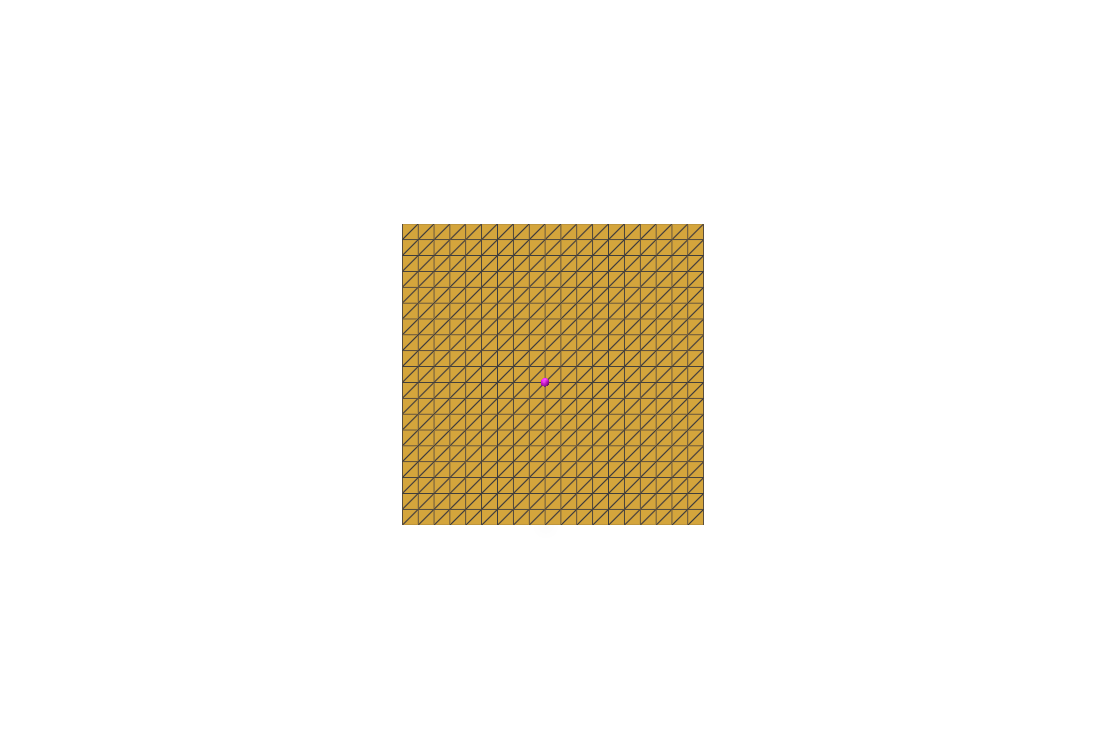}%
		\includegraphics[trim={12cm 5cm 12cm 5cm},clip,width=0.32\textwidth]{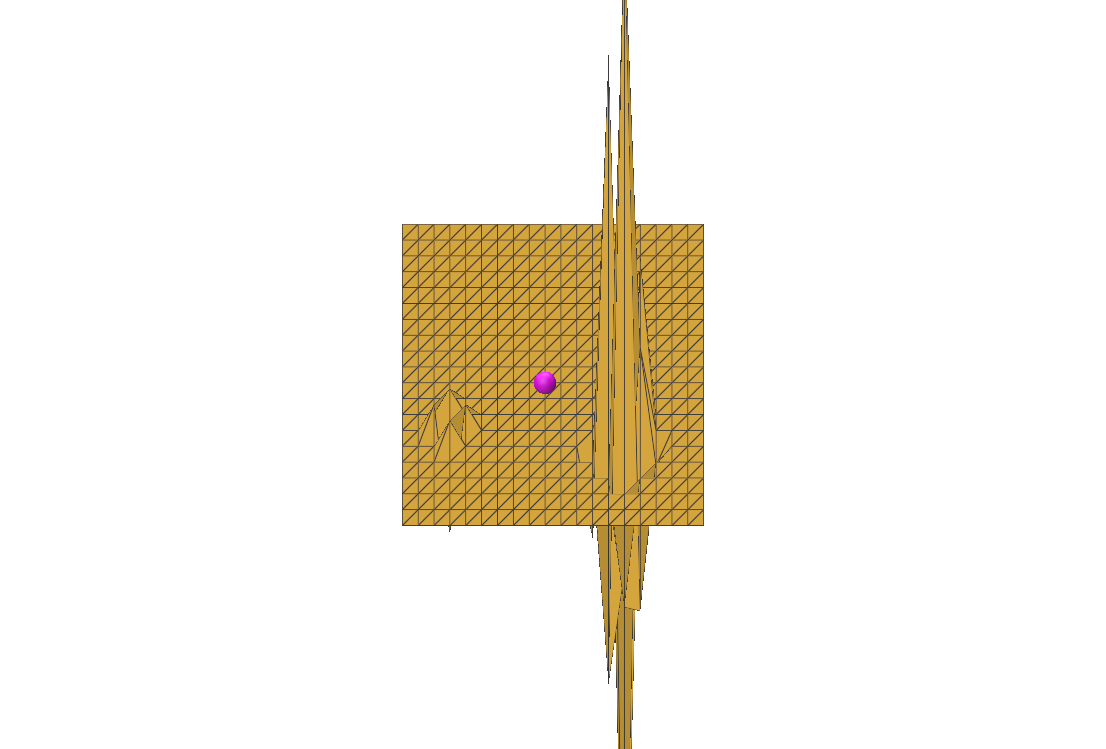}%
		\blankgraphics{2\wNine}{2.8\wNine}%
		%\hspace{4pt} 
\end{minipage}

		\caption{DEIM with POD basis}
		\end{subfigure}

		\end{minipage}
		\end{adjustbox}

		\caption{Testing reduced bending projection at different reduced dimensions (6, 10, 20) against full order model. Experiments (left) free falling under gravitational force, and (right) vertex poking under no gravity.}
		\label{fig:bending_only_sim_comparision}
	\end{figure*}

\begin{figure*}[h!]

		\centering
		% tweak horizontal gaps between panels (column padding)
		\setlength{\tabcolsep}{3pt}
		% remove extra row height in the tabular
		\renewcommand{\arraystretch}{0}

		% inside a figure or figure* environment
		\centering
		\begin{adjustbox}{width=\textwidth} % or =\textwidth, =0.8\textwidth, etc.
		\begin{minipage}{\textwidth}

			% ---------- Row 1 ----------
		\begin{subfigure}[t]{\textwidth}
		\centering
		\includegraphics[trim={10cm 7cm 10cm 7cm},clip,width=0.4\textwidth]{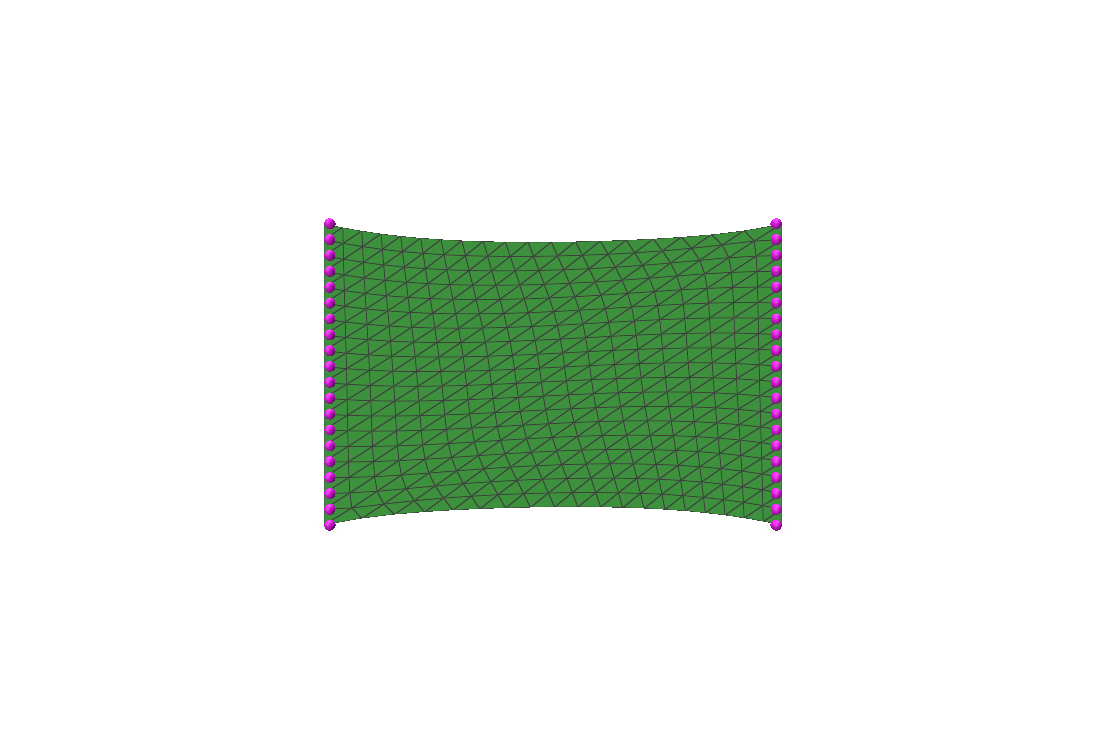}\hfill
		\includegraphics[trim={12cm 9cm 12cm 7cm},clip,width=0.4\textwidth]{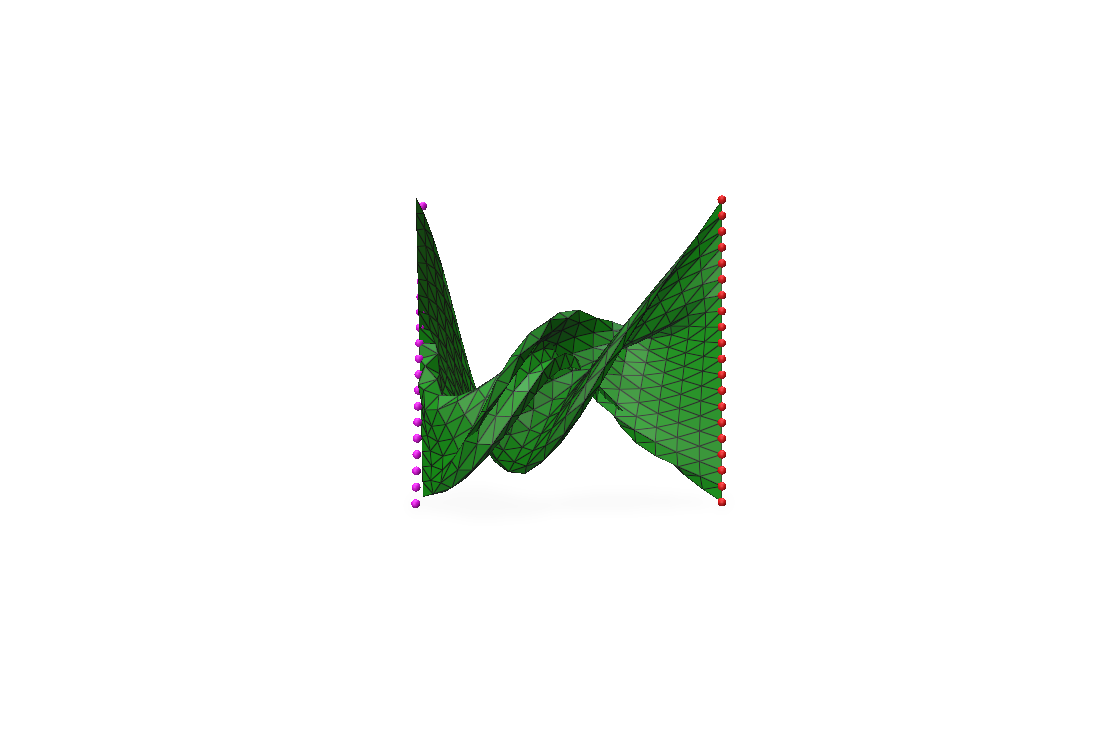}
		
		\caption{FOM}
		\end{subfigure}
		\vspace{-4pt}
		% ---------- Row 2 ----------
		\begin{subfigure}[t]{\textwidth}
		\centering
		\newcommand{\wNine}{\dimexpr\linewidth/9\relax}

		\begin{minipage}{0.45\textwidth}
			\centering
			\includegraphics[trim={8cm 7cm 12cm 7cm},clip,width=0.31\linewidth]{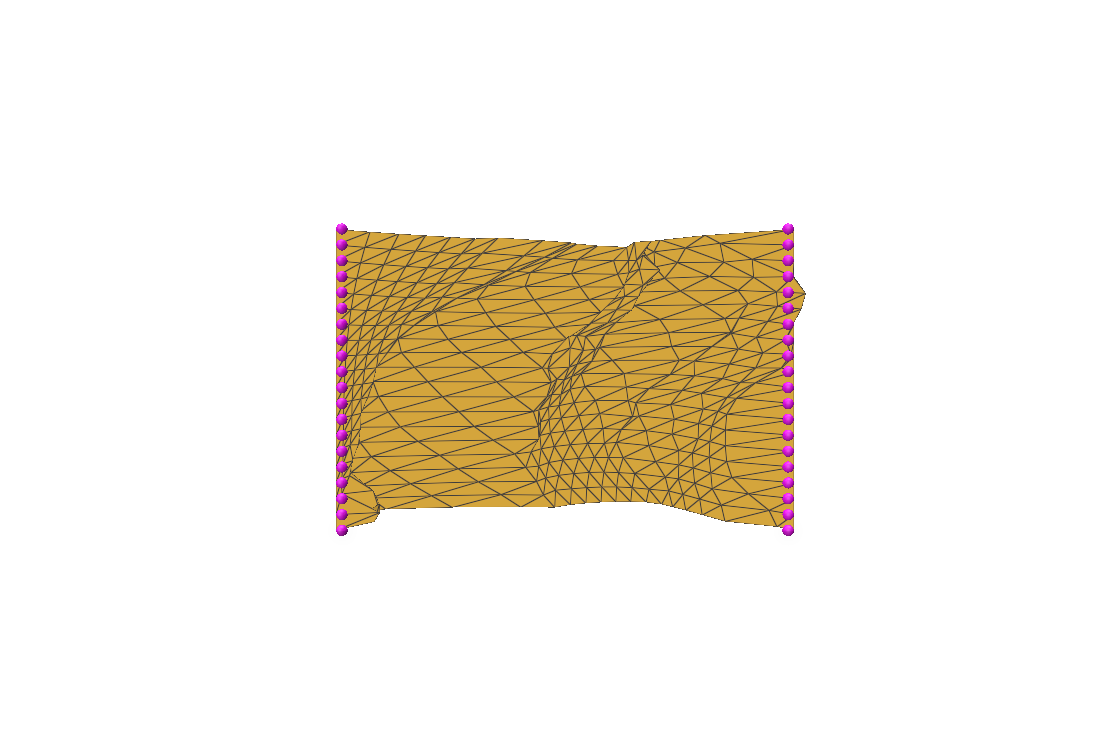}\hfill
			\includegraphics[trim={12cm 7cm 8cm 7cm},clip,width=0.31\linewidth]{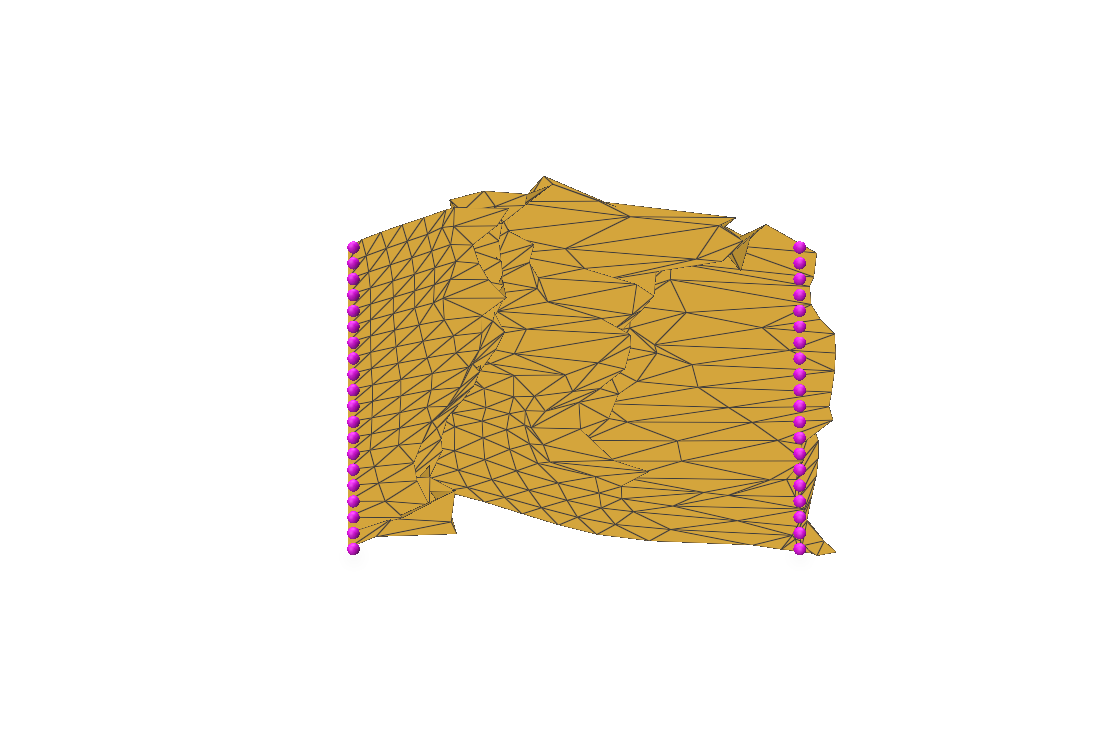}\hfill
			\includegraphics[trim={7cm 5cm 4cm 7cm},clip,width=0.36\linewidth]{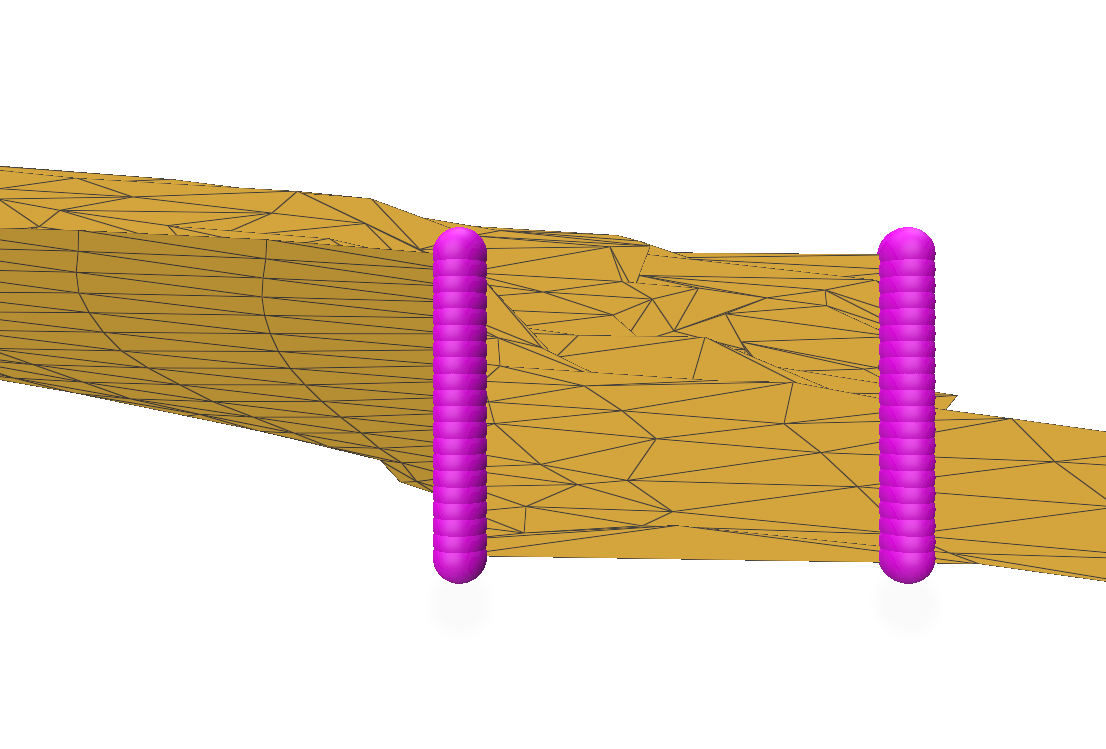}
		\end{minipage}
			\hfill
			% Right triplet
		\begin{minipage}{0.45\textwidth}
			\centering
			\includegraphics[trim={12cm 5cm 12cm 0cm},clip,width=0.31\linewidth]{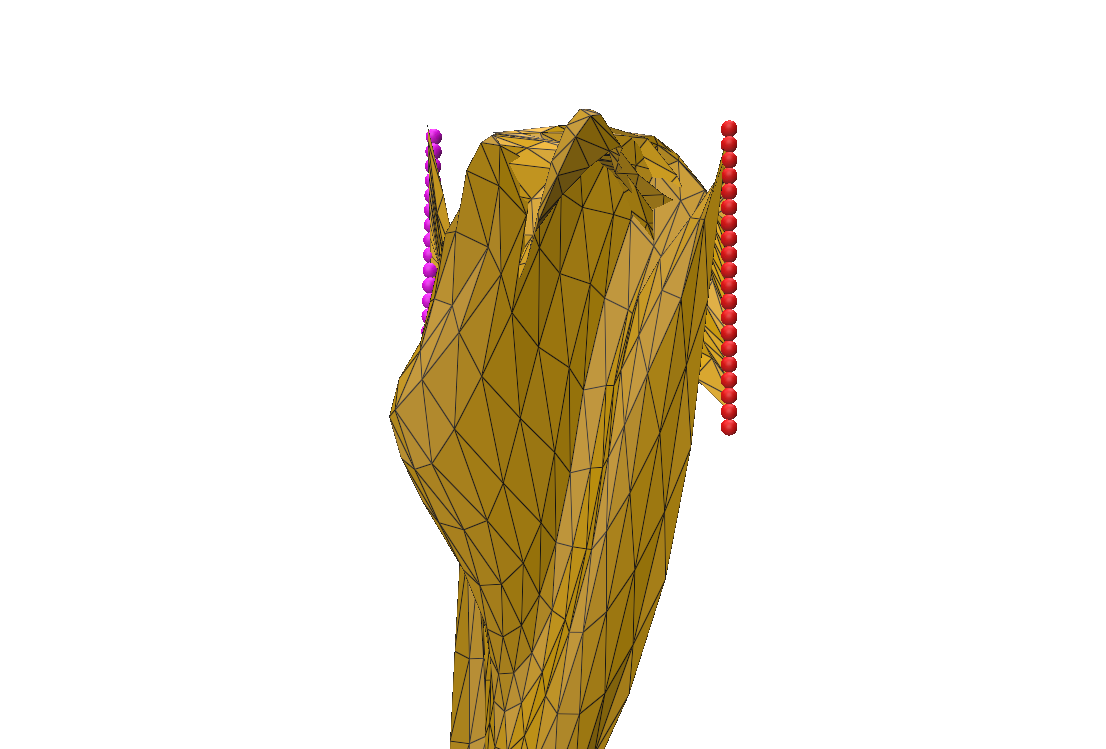}\hfill
			\includegraphics[trim={8cm 5cm 8cm 5cm},clip,width=0.36\linewidth]{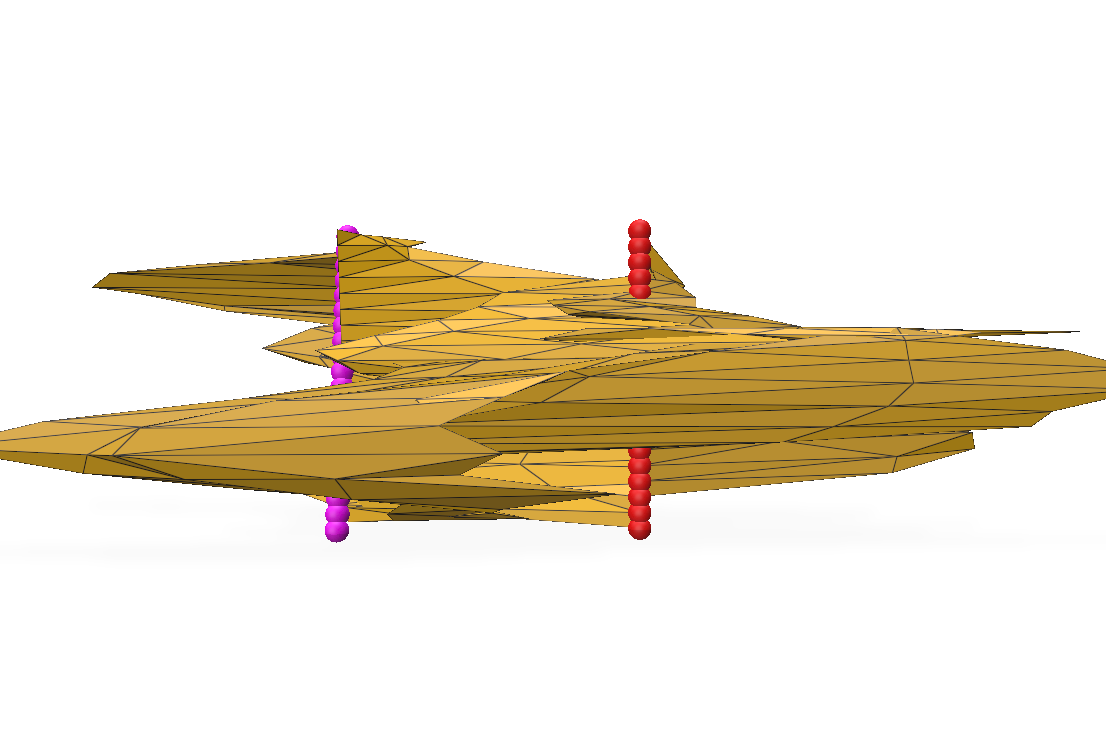}\hfill
			\includegraphics[trim={12cm 5cm 12cm 5cm},clip,width=0.31\linewidth]{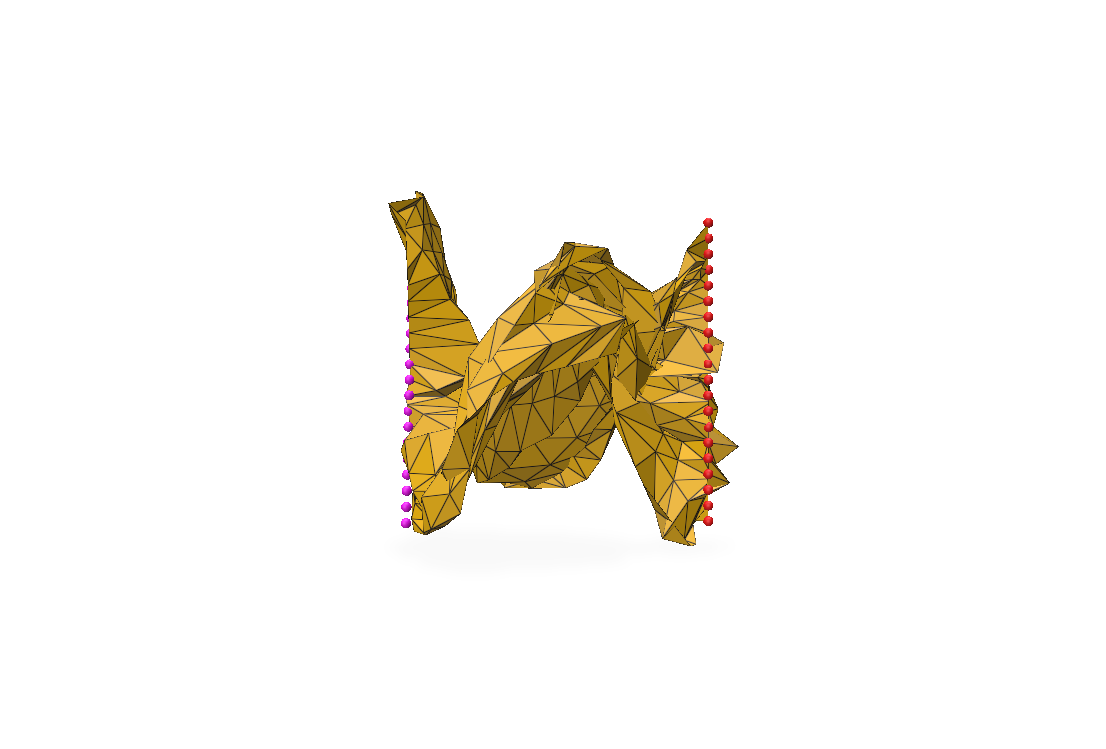}
		\end{minipage}

		\caption{DEIM with POD basis}
		\end{subfigure}

			\end{minipage}
		\end{adjustbox}

		\caption{Testing reduced spring constraint projection at different reduced dimensions (6, 10, 20) against full order model. Experiments (left) stretching, and (right) twisting.}
		\label{fig:spring_only_sim_comparision}
	\end{figure*}

\begin{figure*}[h!]

		\centering
		% tweak horizontal gaps between panels (column padding)
		\setlength{\tabcolsep}{3pt}
		% remove extra row height in the tabular
		\renewcommand{\arraystretch}{0}

		% inside a figure or figure* environment
		\centering
		\begin{adjustbox}{width=\textwidth} % or =\textwidth, =0.8\textwidth, etc.
		\begin{minipage}{\textwidth}

			% ---------- Row 1 ----------
		\begin{subfigure}[t]{\textwidth}
		\centering
		\includegraphics[trim={10cm 7cm 10cm 7cm},clip,width=0.4\textwidth]{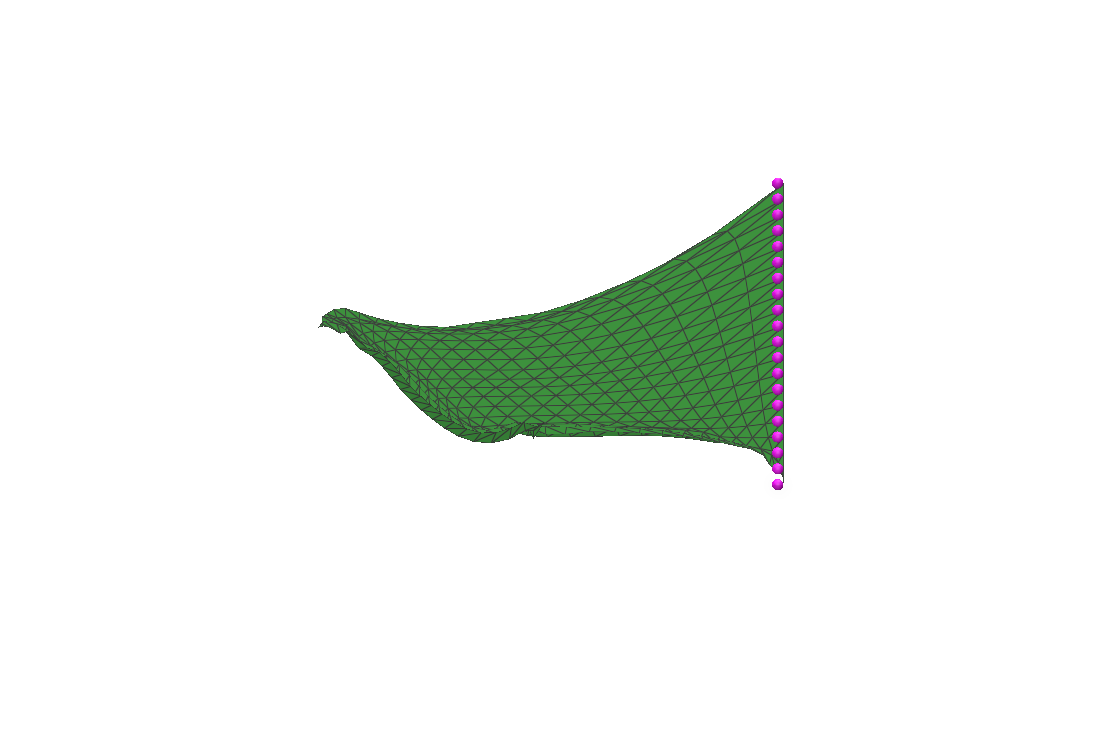}\hfill
		\includegraphics[trim={12cm 9cm 12cm 7cm},clip,width=0.4\textwidth]{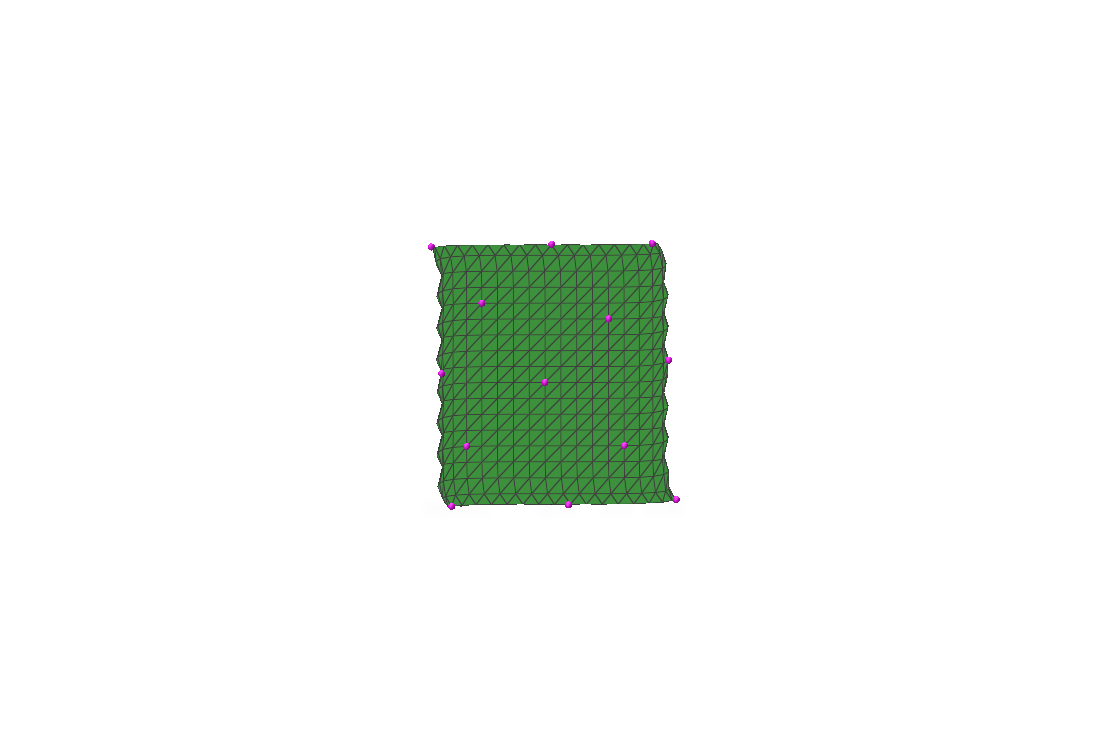}
		
		\caption{FOM}
		\end{subfigure}
		\vspace{-4pt}
		% ---------- Row 2 ----------
		\begin{subfigure}[t]{\textwidth}
		\centering
		\newcommand{\wNine}{\dimexpr\linewidth/9\relax}

		\begin{minipage}{0.45\textwidth}
			\centering
			\includegraphics[trim={8cm 5cm 11cm 5cm},clip,width=0.32\linewidth]{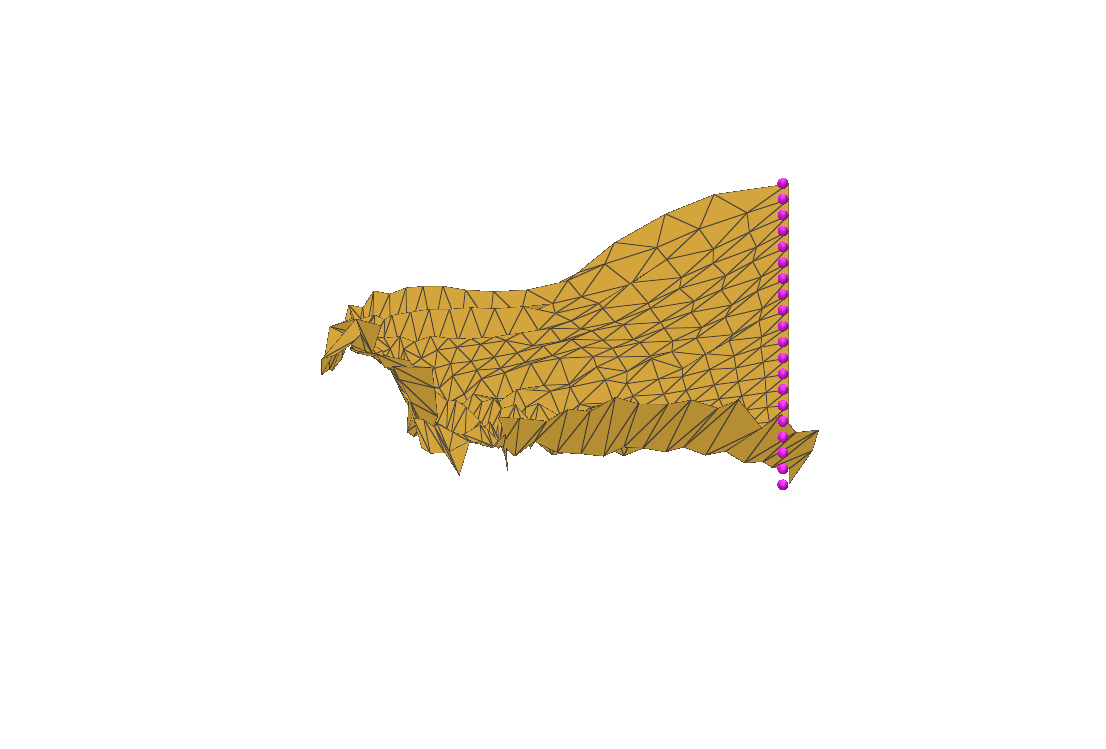}\hfill
			\includegraphics[trim={10cm 5cm 8cm 5cm},clip,width=0.32\linewidth]{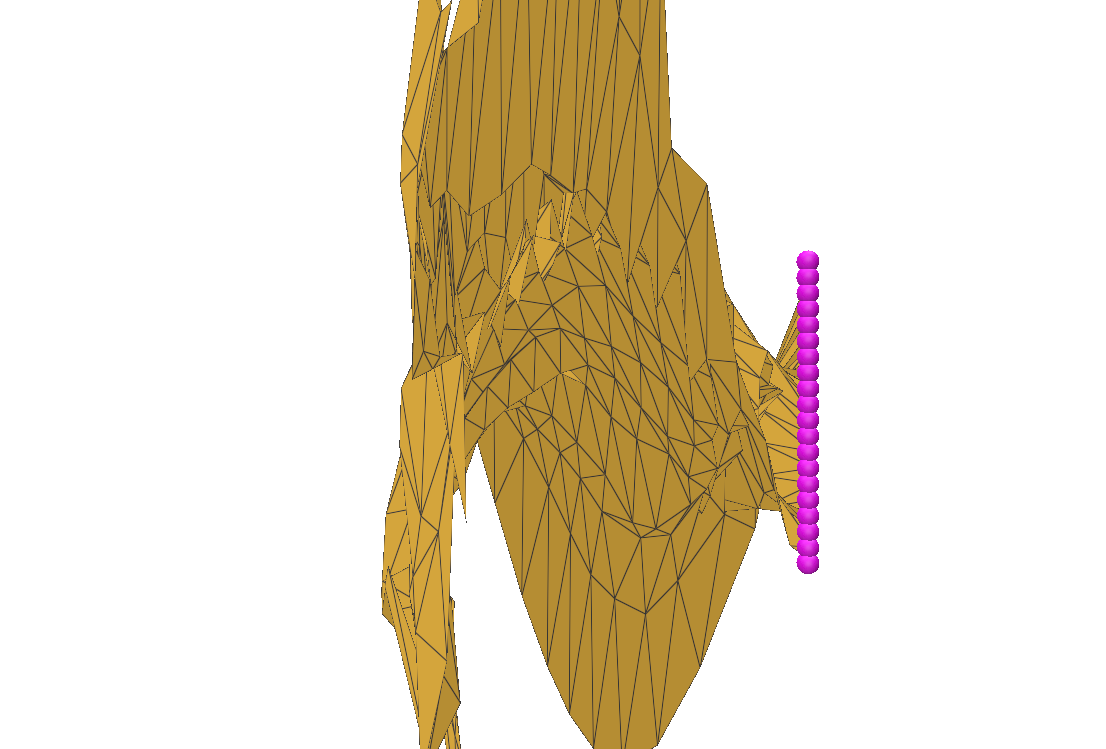}\hfill
			\includegraphics[trim={8cm 5cm 8cm 5cm},clip,width=0.32\linewidth]{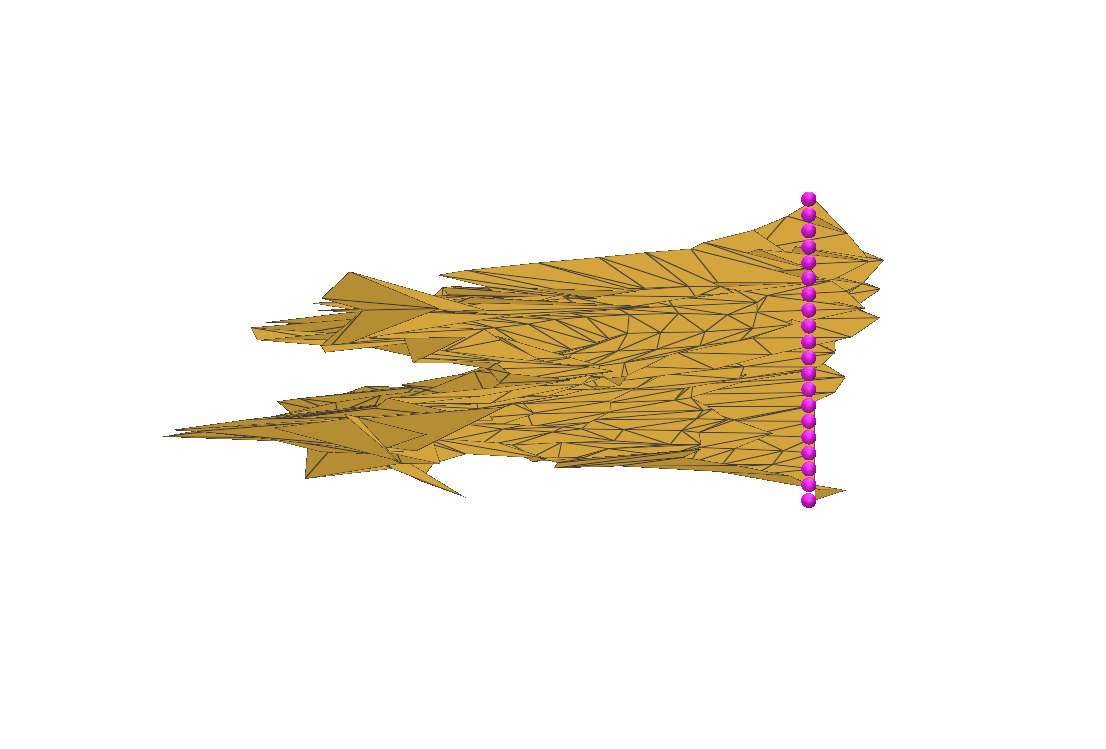}
		\end{minipage}
			\hfill
			% Right triplet
		\begin{minipage}{0.45\textwidth}
			\centering
			\includegraphics[trim={12cm 5cm 12cm 0cm},clip,width=0.32\linewidth]{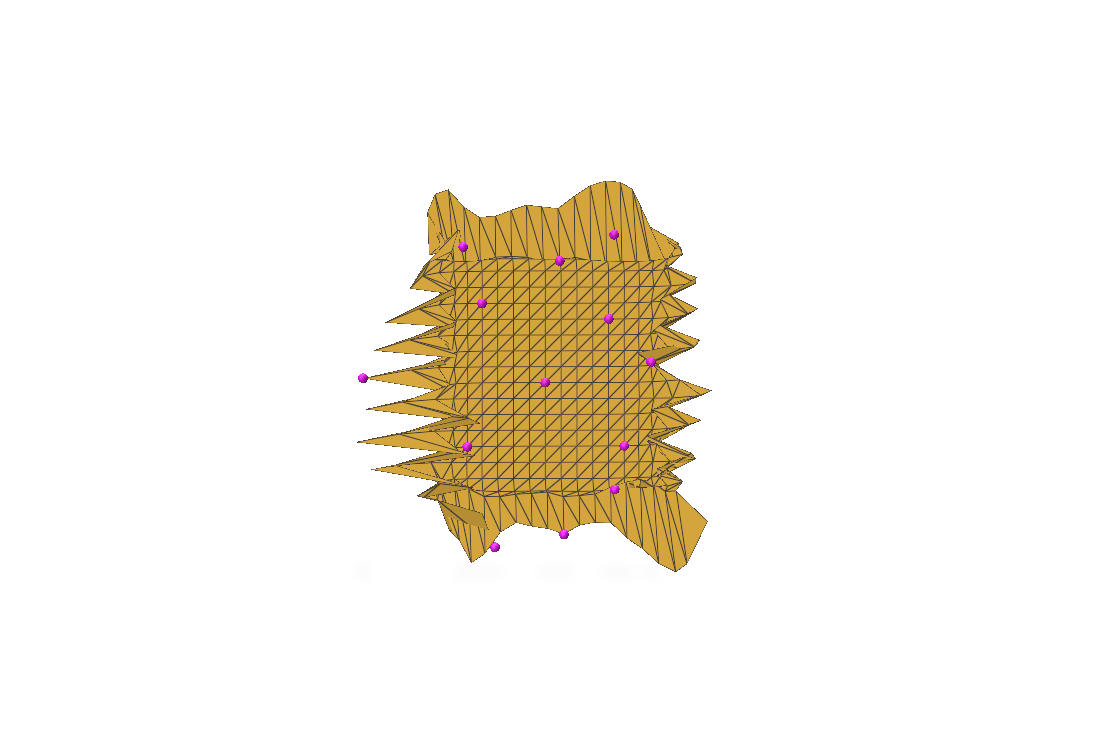}\hfill
			\includegraphics[trim={12cm 5cm 12cm 0cm},clip,width=0.32\linewidth]{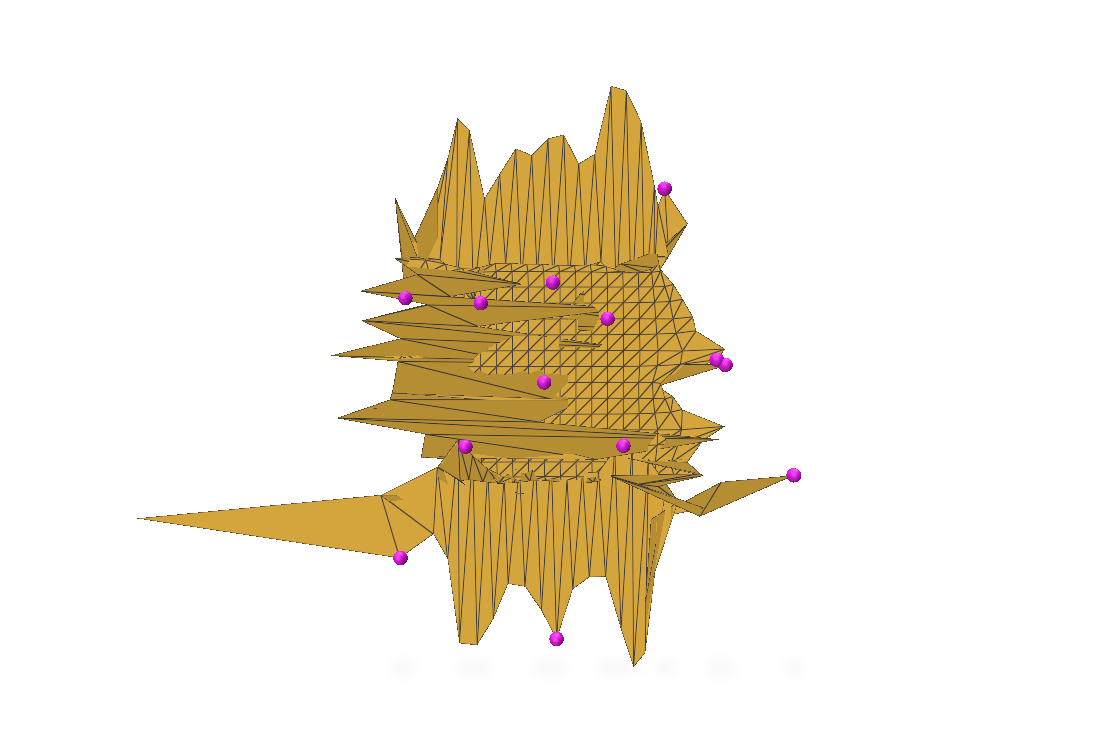}\hfill
			\includegraphics[trim={12cm 5cm 12cm 0cm},clip,width=0.32\linewidth]{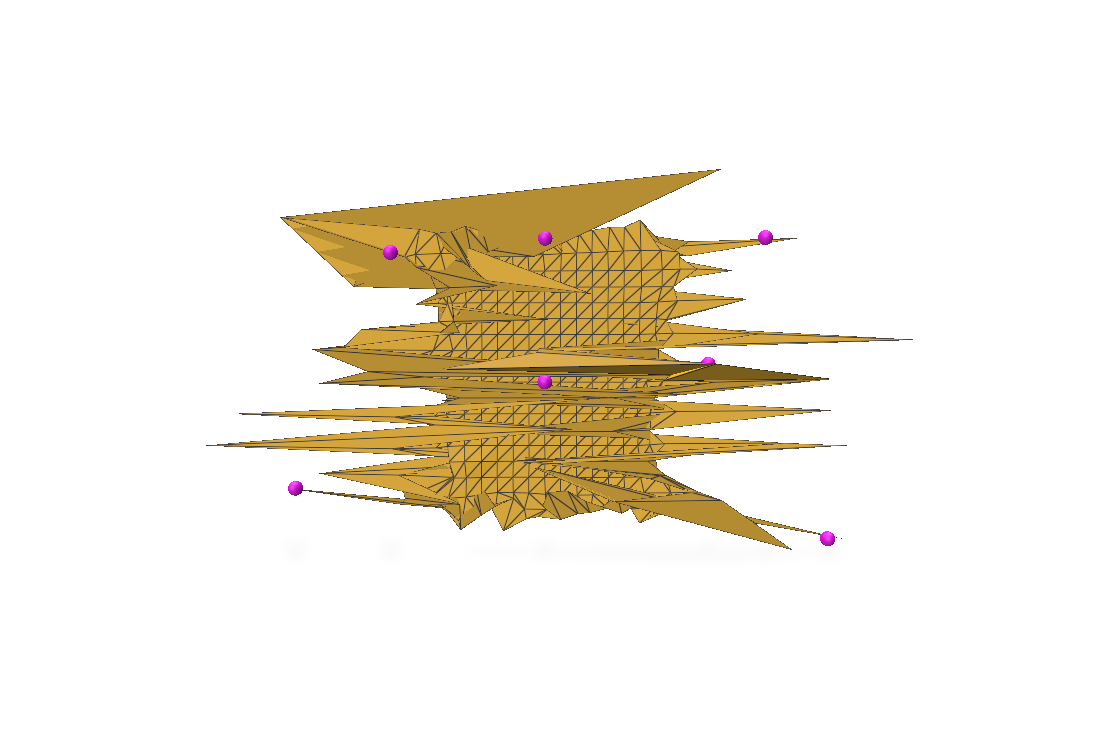}
		\end{minipage}

		\caption{DEIM with POD basis}
		\end{subfigure}

			\end{minipage}
		\end{adjustbox}

		\caption{Testing reduced triangle strain constraint projection at different reduced dimensions (6, 10, 20)against full order model. Experiments (left) stretching, and (right) poking.}
		\label{fig:strain_only_sim_comparision}
	\end{figure*}

	This behavior prompted a deeper investigation. It was unclear why the reduced solver performed reasonably well for bending (defined around vertices) and triangle strain (defined on faces), but failed for edge springs. To clarify this, we designed separate, minimal examples for each constraint type, i.e., simulations in which only one type of constraint was active (see Figures~\ref{fig:bending_only_sim_comparision},~\ref{fig:spring_only_sim_comparision}, ~\ref{fig:strain_only_sim_comparision}).

	Designing a meaningful test for bending alone, without enforcing edge-length constraints, was particularly challenging (see Figure~\ref{fig:bending_only_sim_comparision}), we tried simple gravitational free fall and vertex poking. As a result, we constructed a naive test where a cloth sheet falls vertically under gravity, constrained only by bending. Even in this simple setting, the reduced solver proved highly unstable. As we increased the number of basis vectors, the simulation failed numerically for most values of reduced dimension, and succeeded only irregularly.

	For triangle strain, Figure~\ref{fig:strain_only_sim_comparision}, the full (non-reduced) simulation behaved plausibly in both tests of stretching and twisting under strain tension. However, the reduced simulation exhibited instability and strong visual artifacts, even when using the full number of basis vectors computed via SVD, without truncation. Similarly, reducing the spring forces led to severely distorted simulations (see Figure~\ref{fig:spring_only_sim_comparision}), and artifacts grow with the number of basis.

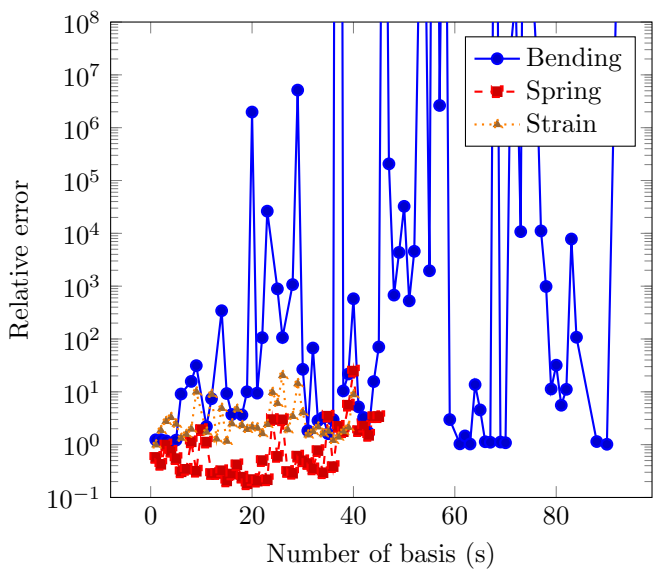
\begin{figure}[t!]		
			\begin{tikzpicture}
				\begin{axis}[
					width=\linewidth,            % fits the column
					height=0.9\linewidth,       % nice aspect for columns
					xlabel={Number of basis (s)},
					ylabel={Relative error},
					ymode=log,               % make Y axis logarithmic
					ymin=1e-1, ymax=1e8,     % adjust to data range
					% % grid=both,
					% minor grid style={gray!15},
					% major grid style={gray!35},
					legend cell align=left,
					legend pos=north east,
					% If your CSVs use semicolons, add: table/col sep=semicolon,
					]

					% --- Series 1: data1.csv with columns 'x' and 'y' ---
					\addplot+[y filter/.expression={and(\thisrow{relative_errors_z}>=1e-2, \thisrow{relative_errors_z}<=1e2) ? y : nan},
					thick, mark=*,
					] table[
					x=numPoints, y=relative_errors_z,
					col sep=comma,
					] {cloth_verts_bending_deim_pod_vectorized_train_convergence_tests.csv};
					\addlegendentry{Bending}

					% --- Series 2: data2.csv with columns 'x' and 'y' ---
					\addplot+[
					thick, dashed, mark=square*,
					] table[
					x=numPoints, y=relative_errors_x,
					col sep=comma,
					] {cloth_edge_spring_deim_pod_vectorized_train_convergence_tests.csv};
					\addlegendentry{Spring}

					% % --- Series 3: data3.csv with columns 'x' and 'y' ---
					\addplot+[
					thick, dotted, mark=triangle*, color=orange
					] table[
					x=numPoints, y=relative_errors_x,
					col sep=comma,
					] {cloth_tris_strain_deim_pod_vectorized_train_convergence_tests.csv};
					\addlegendentry{Strain}

					% % --- Series 1: data1.csv with columns 'x' and 'y' ---
					% \addplot+[
					% thick, color=orange
					% ] table[
					% x=numPoints, y=relative_errors_x,
					% col sep=comma,
					% ] {bar_tets_deformation_gradient_deim_pod_train_convergence_tests.csv};
					% \addlegendentry{Deformation gradient}

				\end{axis}
			\end{tikzpicture}
			\caption{Relative reconstruction error $\frac{\norm{p(t)- V \hat{p}(t)}}{\norm{p(t)}}$ measured for different reduced constraint projections for POD basis combined with DEIM, related to simulations shown in figure (\ref{fig:deim_pod_singVals}).
	.}
			\label{fig:deim_reconstruction_relative_error}
	\end{figure}

\begin{figure}[t!]		
			\begin{tikzpicture}
				\begin{axis}[
					width=\linewidth,            % fits the column
					height=0.9\linewidth,       % nice aspect for columns
					xlabel={Number of basis (s)},
					ylabel={Relative error},
					ymode=log,               % make Y axis logarithmic
					ymin=0, ymax=500,     % adjust to data range
					% % grid=both,
					% minor grid style={gray!15},
					% major grid style={gray!35},
					legend cell align=left,
					legend pos=north east % If your CSVs use semicolons, add: table/col sep=semicolon,
					]

					% --- Series 1: data1.csv with columns 'x' and 'y' ---
					\addplot+[
					thick, mark=*,
					] table[
					x=component, y=singVal,
					col sep=comma,
					] {cloth_verts_bending_constrprojBases_pcaExtraction_singValues.csv};
					\addlegendentry{Bending-poking}

					% --- Series 2: data2.csv with columns 'x' and 'y' ---
					\addplot+[
					thick, dashed, mark=square*,
					] table[
					x=component, y=singVal,
					col sep=comma,
					] {cloth_edge_spring_constrprojBases_pcaExtraction_singValues.csv};
					\addlegendentry{Spring-twisting}

					% % --- Series 3: data3.csv with columns 'x' and 'y' ---
					\addplot+[
					thick, dotted, mark=triangle*, color=orange
					] table[
					x=component, y=singVal,
					col sep=comma,
					] {cloth_tris_strain_constrprojBases_pcaExtraction_singValues.csv};
					\addlegendentry{Strain-stretching}

					% % --- Series 1: data1.csv with columns 'x' and 'y' ---
					% \addplot+[
					% thick, color=orange
					% ] table[
					% x=numPoints, y=relative_errors_x,
					% col sep=comma,
					% ] {bar_tets_deformation_gradient_deim_pod_train_convergence_tests.csv};
					% \addlegendentry{Deformation gradient}

				\end{axis}
			\end{tikzpicture}
			\caption{Singular values computed for some of the experiments shown in Figures (\ref{fig:bending_only_sim_comparision}, \ref{fig:spring_only_sim_comparision}, \ref{fig:strain_only_sim_comparision}) }

			\label{fig:deim_pod_singVals}
	\end{figure}
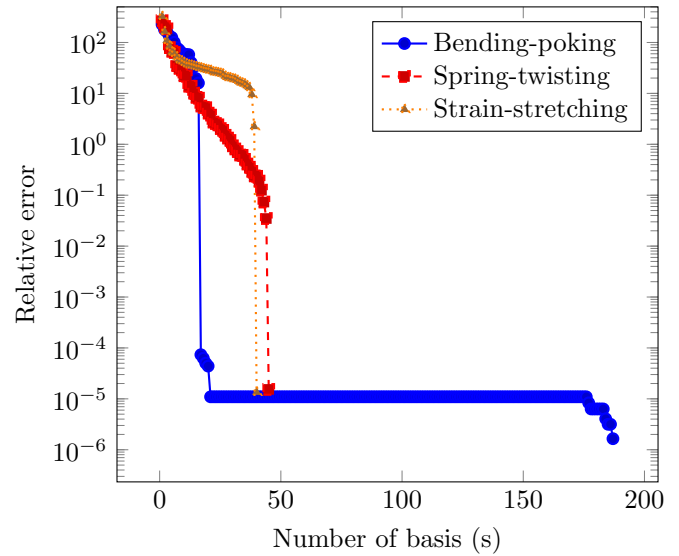

	These findings suggest that, in the first experiment shown in Figure~\ref{fig:cloth_sim_grid}, the spring constraints on edges were the dominant force governing the dynamics. When spring forces are solved in the full space, reduced bending or strain bases may still introduce tolerable visual artifacts. However, once spring forces are also reduced, the simulation becomes unstable and breaks down completely.
	The POD basis, combined with DEIM interpolation, failed to produce representative simulations—even for the same dynamics on which they were trained. Although in figure ~\ref{fig:deim_pod_singVals} one observe nicely decaying singular value associated to the vectorized POD basis, as shown in Figure~\ref{fig:deim_reconstruction_relative_error}, the reconstruction error on the training snapshots remains high.  This suggests that simply computing more basis vectors from additional snapshots may not resolve the issue.

	An important consideration in this context is that constraint projections are computed pointwise, which makes them highly amenable to parallelization. As more basis vectors are included, the DEIM solver matrix \( \left( P^\top V \right) \) becomes denser, and at some point, the computational cost of solving this system outweighs the benefit, particularly when compared to applying parallelized computations across all elements.

% ----------------------------------------------------------	

	% All constraints are in volve in the simulations, but they are reduced one by one separately
% ----------------------------------------------------------	
	% Comparison between FOM and reduced for simulations with one constraint type

% ----------------------------------------------------------	

	% reconstruction errors for deim pod for separately simulated constraints

	% POD singularvals for deim pod for separately simulated constraints

%%%%%%%%%%%%%%%%%%%%%%%%%%%%%%%%%%%%%%%%%%%%%%%%%%%%%%%%%%%%%%%%%%%%%%%%%%%%%%%%
%%%%%%%%%%%%%%%%%%%%%%%%%%%%%%%%%%%%%%%%%%%%%%%%%%%%%%%%%%%%%%%%%%%%%%%%%%%%%%%%
\section{Conclusion}

	Establishing that a method works reliably in all cases requires rigorous proof, while demonstrating its limitations may take only a single counter example at least in its current form of development. In this study, we evaluated the compatibility of POD-computed bases and the capability of the DEIM algorithm to select a small subset of elements for computing nonlinear forces in real time.

	Our main observation is that POD bases failed to capture suitable reduced subspaces for constraint manifolds. This was demonstrated in the second experiment, where each constraint type was simulated individually, time-varying snapshots were collected, and reduced bases were computed separately. In each case, the resulting reduced simulations exhibited either significant numerical instability, visual artifacts, or complete solver failure.

	One might still attempt to reduce non-dominant constraint forces while solving dominant constraints in full dimension, as in the first experiment. However, this approach introduces non-intuitive design choices and often results in large, dense solvers for the reduced portion, potentially outweighing the benefits compared to parallel evaluation of constraint projections on all elements.

	In summary, based on our experiments, using DEIM to accelerate nonlinear constraint projections for physics-based animation is not recommended. By presenting these results, we aim to highlight the current limitations of the method in this context, with the hope that it may inspire further development and refinement within the community.

%%%%%%%%%%%%%%%%%%%%%%%%%%%%%%%%%%%%%%%%%%%%%%%%%%%%%%%%%%%%%%%%%%%%%%%%%%%%%%%%
%%%%%%%%%%%%%%%%%%%%%%%%%%%%%%%%%%%%%%%%%%%%%%%%%%%%%%%%%%%%%%%%%%%%%%%%%%%%%%%%

\section{Acknowledgement}%
\label{sec:ack}                                                      Authors would like to acknowledge the German Research Foundation (DFG) Research Training Group 2297 ”MathCoRe”, Magdeburg.
%%%%%%%%%%%%%%%%%%%%%%%%%%%%%%%%%%%%%%%%%%%%%%%%%%%%%%%%%%%%%%%%%%%%%%%%%%%%%%%%

\addcontentsline{toc}{section}{References}
\bibliographystyle{plainurl}
\bibliography{exampleref}

@article{Brandt18,
  title={Hyper-reduced projective dynamics},
  author={Brandt, Christopher and Eisemann, Elmar and Hildebrandt, Klaus},
  journal={ACM Transactions on Graphics (TOG)},
  volume={37},
  number={4},
  pages={1--13},
  year={2018},
  publisher={ACM New York, NY, USA}
}

@article{Bouaziz14,
  title={Projective dynamics: Fusing constraint projections for fast simulation},
  author={Bouaziz, Sofien and Martin, Sebastian and Liu, Tiantian and Kavan, Ladislav and Pauly, Mark},
  journal={ACM Transactions on graphics (TOG)},
  volume={33},
  number={4},
  pages={1--11},
  year={2014},
  publisher={ACM New York, NY, USA}
}

@misc{DiffPD21,
      title={DiffPD: Differentiable Projective Dynamics}, 
      author={Tao Du and Kui Wu and Pingchuan Ma and Sebastien Wah and Andrew Spielberg and Daniela Rus and Wojciech Matusik},
      year={2021},
      eprint={2101.05917},
      archivePrefix={arXiv},
      primaryClass={cs.LG},
      url={https://arxiv.org/abs/2101.05917}, 
}

@ARTICLE{Xiaowei18,
  author={He, Xiaowei and Wang, Huamin and Wu, Enhua},
  journal={IEEE Transactions on Visualization and Computer Graphics}, 
  title={Projective Peridynamics for Modeling Versatile Elastoplastic Materials}, 
  year={2018},
  volume={24},
  number={9},
  pages={2589-2599},
  doi={10.1109/TVCG.2017.2755646}}

@article{Zhendong18,
author = {Wang, Zhendong and Wu, Longhua and Fratarcangeli, Marco and Tang, Min and Wang, Huamin},
title = {Parallel Multigrid for Nonlinear Cloth Simulation},
journal = {Computer Graphics Forum},
volume = {37},
number = {7},
pages = {131-141},
doi = {https://doi.org/10.1111/cgf.13554},
url = {https://onlinelibrary.wiley.com/doi/abs/10.1111/cgf.13554},
eprint = {https://onlinelibrary.wiley.com/doi/pdf/10.1111/cgf.13554},
year = {2018}
}

@article{LTB17,
author = {Liu, Tiantian and Bouaziz, Sofien and Kavan, Ladislav},
title = {Quasi Newton Methods for Real-Time Simulation of Hyperelastic Materials},
year = {2017},
issue_date = {August 2017},
publisher = {Association for Computing Machinery},
address = {New York, NY, USA},
volume = {36},
number = {4},
issn = {0730-0301},
url = {https://doi.org/10.1145/3072959.2990496},
doi = {10.1145/3072959.2990496},
journal = {ACM Trans. Graph.},
month = jul,
articleno = {116a},
numpages = {16}
}

@article{huang2019survey,
  title={A survey on fast simulation of elastic objects},
  author={Huang, Jin and Chen, Jiong and Xu, Weiwei and Bao, Hujun},
  journal={Frontiers of Computer Science},
  volume={13},
  pages={443--459},
  year={2019},
  publisher={Springer}
}

@article{Pentland_Williams89,
author = {Pentland, A. and Williams, J.},
title = {Good Vibrations: Modal Dynamics for Graphics and Animation},
year = {1989},
issue_date = {July 1989},
publisher = {Association for Computing Machinery},
address = {New York, NY, USA},
volume = {23},
number = {3},
issn = {0097-8930},
url = {https://doi.org/10.1145/74334.74355},
doi = {10.1145/74334.74355},
journal = {SIGGRAPH Comput. Graph.},
month = {jul},
pages = {207–214},
numpages = {8}
}

@inproceedings{PBS15,
  title={Position-Based Simulation Methods in Computer Graphics.},
  author={Bender, Jan and M{\"u}ller, Matthias and Macklin, Miles},
  booktitle={Eurographics (tutorials)},
  pages={8},
  year={2015}
}

@inproceedings{MMMCN16,
author = {Macklin, Miles and M\"{u}ller, Matthias and Chentanez, Nuttapong},
title = {XPBD: position-based simulation of compliant constrained dynamics},
year = {2016},
isbn = {9781450345927},
publisher = {Association for Computing Machinery},
address = {New York, NY, USA},
url = {https://doi.org/10.1145/2994258.2994272},
doi = {10.1145/2994258.2994272},
booktitle = {Proceedings of the 9th International Conference on Motion in Games},
pages = {49–54},
numpages = {6},
location = {Burlingame, California},
series = {MIG '16}
}

@article{DKWB2018,
author = {Crispin Deul and Tassilo Kugelstadt and Marcel Weiler and Jan Bender},
title = {Direct Position-Based Solver for Stiff Rods},
year = {2018},
journal = {Computer Graphics Forum},
volume = {37},
number = {6},
pages = {313-324},
doi = {10.1111/cgf.13326},
url = {https://onlinelibrary.wiley.com/doi/abs/10.1111/cgf.13326},
eprint = {https://onlinelibrary.wiley.com/doi/pdf/10.1111/cgf.13326},
}

@inproceedings{BDWA98,
author = {Baraff, David and Witkin, Andrew},
title = {Large steps in cloth simulation},
year = {1998},
isbn = {0897919998},
publisher = {Association for Computing Machinery},
address = {New York, NY, USA},
url = {https://doi.org/10.1145/280814.280821},
doi = {10.1145/280814.280821},
booktitle = {Proceedings of the 25th Annual Conference on Computer Graphics and Interactive Techniques},
pages = {43–54},
numpages = {12},
series = {SIGGRAPH '98}
}

@article{Zhang19,
  title={Deformable Models for Surgical Simulation: A Survey},
  author={Jinao Zhang and Yongmin Zhong and Chengfan Gu},
  journal={IEEE Reviews in Biomedical Engineering},
  year={2019},
  volume={11},
  pages={143-164},
  url={https://api.semanticscholar.org/CorpusID:50785173}
}

@inproceedings{PBSSurvay14,
  title={A survey on position-based simulation methods in computer graphics},
  author={Bender, Jan and M{\"u}ller, Matthias and Otaduy, Miguel A and Teschner, Matthias and Macklin, Miles},
  booktitle={Computer graphics forum},
  volume={33},
  number={6},
  pages={228--251},
  year={2014},
  organization={Wiley Online Library}
}

@article{snowden2017methods,
  title={Methods of model reduction for large-scale biological systems: a survey of current methods and trends},
  author={Snowden, Thomas J and van der Graaf, Piet H and Tindall, Marcus J},
  journal={Bulletin of mathematical biology},
  volume={79},
  pages={1449--1486},
  year={2017},
  publisher={Springer}
}

@article{Barbic05,
author = {Barbi\v{c}, Jernej and James, Doug L.},
title = {Real-Time Subspace Integration for St. Venant-Kirchhoff Deformable Models},
year = {2005},
issue_date = {July 2005},
publisher = {Association for Computing Machinery},
address = {New York, NY, USA},
volume = {24},
number = {3},
issn = {0730-0301},
url = {https://doi.org/10.1145/1073204.1073300},
doi = {10.1145/1073204.1073300},
journal = {ACM Trans. Graph.},
month = {jul},
pages = {982–990},
numpages = {9}
}

@inproceedings{Terzopoulos87,
author = {Terzopoulos, Demetri and Platt, John and Barr, Alan and Fleischer, Kurt},
title = {Elastically deformable models},
year = {1987},
isbn = {0897912276},
publisher = {Association for Computing Machinery},
address = {New York, NY, USA},
url = {https://doi.org/10.1145/37401.37427},
doi = {10.1145/37401.37427},
booktitle = {Proceedings of the 14th Annual Conference on Computer Graphics and Interactive Techniques},
pages = {205–214},
numpages = {10},
series = {SIGGRAPH '87}
}

@inproceedings{Trusty23,
author = {Trusty, Ty and Benchekroun, Otman and Grinspun, Eitan and Kaufman, Danny M. and Levin, David I.W.},
title = {Subspace Mixed Finite Elements for Real-Time Heterogeneous Elastodynamics},
year = {2023},
isbn = {9798400703157},
publisher = {Association for Computing Machinery},
address = {New York, NY, USA},
url = {https://doi.org/10.1145/3610548.3618220},
doi = {10.1145/3610548.3618220},
booktitle = {SIGGRAPH Asia 2023 Conference Papers},
articleno = {112},
numpages = {10},
location = {<conf-loc>, <city>Sydney</city>, <state>NSW</state>, <country>Australia</country>, </conf-loc>},
series = {SA '23}
}

@article{Brandt19,
author = {Brandt, Christopher and Scandolo, Leonardo and Eisemann, Elmar and Hildebrandt, Klaus},
title = {The reduced immersed method for real-time fluid-elastic solid interaction and contact simulation},
year = {2019},
issue_date = {December 2019},
publisher = {Association for Computing Machinery},
address = {New York, NY, USA},
volume = {38},
number = {6},
issn = {0730-0301},
url = {https://doi.org/10.1145/3355089.3356496},
doi = {10.1145/3355089.3356496},
journal = {ACM Trans. Graph.},
month = {nov},
articleno = {191},
numpages = {16}
}

@inproceedings{Barbivc11,
author = {Barbi\v{c}, Jernej and Zhao, Yili},
title = {Real-time large-deformation substructuring},
year = {2011},
isbn = {9781450309431},
publisher = {Association for Computing Machinery},
address = {New York, NY, USA},
url = {https://doi.org/10.1145/1964921.1964986},
doi = {10.1145/1964921.1964986},
booktitle = {ACM SIGGRAPH 2011 Papers},
articleno = {91},
numpages = {8},
location = {Vancouver, British Columbia, Canada},
series = {SIGGRAPH '11}
}

@inproceedings{Monem23,
author = {Monem Abdelhafez, Shaimaa and Benner, Peter and Lessig, Christian},
title = {Improved Projective Dynamics Global Using Snapshots-based Reduced Bases},
year = {2023},
isbn = {9798400701528},
publisher = {Association for Computing Machinery},
address = {New York, NY, USA},
url = {https://doi.org/10.1145/3588028.3603665},
doi = {10.1145/3588028.3603665},
booktitle = {ACM SIGGRAPH 2023 Posters},
articleno = {4},
numpages = {2},
location = {Los Angeles, CA, USA},
series = {SIGGRAPH '23}
}

@misc{Monem25,
      title={SIGGRAPH: G: Improved Projective Dynamics Global Using Snapshots-based Reduced Bases}, 
      author={Shaimaa Monem and Peter Benner and Christian Lessig},
      year={2025},
      eprint={2502.07757},
      archivePrefix={arXiv},
      primaryClass={math.DS},
      url={https://arxiv.org/abs/2502.07757}, 
}

@article{GKIMISIS2025118115,
title = {Non-intrusive reduced-order modeling for dynamical systems with spatially localized features},
journal = {Computer Methods in Applied Mechanics and Engineering},
volume = {444},
pages = {118115},
year = {2025},
issn = {0045-7825},
doi = {https://doi.org/10.1016/j.cma.2025.118115},
url = {https://www.sciencedirect.com/science/article/pii/S0045782525003871},
author = {Leonidas Gkimisis and Nicole Aretz and Marco Tezzele and Thomas Richter and Peter Benner and Karen E. Willcox}
}

@inbook{lassila14,
author = {Lassila, Toni and Manzoni, Andrea and Quarteroni, Alfio and Rozza, Gianluigi},
year = {2014},
month = {01},
pages = {},
title = {Model Order Reduction in Fluid Dynamics: Challenges and Perspectives},
volume = {9},
isbn = {978-3-319-02089-1},
doi = {10.1007/978-3-319-02090-7_9}
}

@article{Radulescu12,
  title={Reduction of dynamical biochemical reactions networks in computational biology},
  author={Radulescu, Ovidiu and Gorban, Alexander N and Zinovyev, Andrei and Noel, Vincent},
  journal={Frontiers in genetics},
  volume={3},
  pages={131},
  year={2012},
  publisher={Frontiers Media SA}
}

@InProceedings{Bnneland19,
  author =	{B{\o}nneland, Frederik Meyer and Jensen, Peter Gj{\o}l and Larsen, Kim G. and Mu\~{n}iz, Marco and Srba, Ji\v{r}{\'\i}},
  title =	{{Partial Order Reduction for Reachability Games}},
  booktitle =	{30th International Conference on Concurrency Theory (CONCUR 2019)},
  pages =	{23:1--23:15},
  series =	{Leibniz International Proceedings in Informatics (LIPIcs)},
  ISBN =	{978-3-95977-121-4},
  ISSN =	{1868-8969},
  year =	{2019},
  volume =	{140},
  editor =	{Fokkink, Wan and van Glabbeek, Rob},
  publisher =	{Schloss Dagstuhl -- Leibniz-Zentrum f{\"u}r Informatik},
  address =	{Dagstuhl, Germany},
  URL =		{https://drops-dev.dagstuhl.de/entities/document/10.4230/LIPIcs.CONCUR.2019.23},
  URN =		{urn:nbn:de:0030-drops-109251},
  doi =		{10.4230/LIPIcs.CONCUR.2019.23}
}

@article{Chaturantabut10,
author = {Chaturantabut, Saifon and Sorensen, Danny C.},
title = {Nonlinear Model Reduction via Discrete Empirical Interpolation},
journal = {SIAM Journal on Scientific Computing},
volume = {32},
number = {5},
pages = {2737-2764},
year = {2010},
doi = {10.1137/090766498},

URL = { 
    
        https://doi.org/10.1137/090766498
    
    

},
eprint = { 
    
        https://doi.org/10.1137/090766498
    
    

}
}

@article{Drmavc16,
author = {Drma\v{c}, Zlatko and Gugercin, Serkan},
title = {A New Selection Operator for the Discrete Empirical Interpolation Method---Improved A Priori Error Bound and Extensions},
journal = {SIAM Journal on Scientific Computing},
volume = {38},
number = {2},
pages = {A631-A648},
year = {2016},
doi = {10.1137/15M1019271},

URL = { 
    
        https://doi.org/10.1137/15M1019271
    
    

},
eprint = { 
    
        https://doi.org/10.1137/15M1019271
    
    

}
}

@article{Benjamin20,
author = {Peherstorfer, Benjamin and Drma\v{c}, Zlatko and Gugercin, Serkan},
title = {Stability of Discrete Empirical Interpolation and Gappy Proper Orthogonal Decomposition with Randomized and Deterministic Sampling Points},
journal = {SIAM Journal on Scientific Computing},
volume = {42},
number = {5},
pages = {A2837-A2864},
year = {2020},
doi = {10.1137/19M1307391},

URL = { 
    
        https://doi.org/10.1137/19M1307391
    
    

},
eprint = { 
    
        https://doi.org/10.1137/19M1307391
    
    

}
}

@article{chaturantabut2012state,

  title={A state space error estimate for nonlinear model reduction},

  author={Chaturantabut, Saifon and Sorensen, Danny C},

  journal={SIAM Journal on Scientific Computing},

  volume={50},

  number={1},

  pages={46--63},

  year={2012},

  publisher={SIAM}

}

@article{bremer2017pod,

  title={POD-DEIM for efficient reduction of a dynamic 2D catalytic reactor model},

  author={Bremer, Jens and Goyal, Pawan and Feng, Lihong and Benner, Peter and Sundmacher, Kai},

  journal={Computers \& Chemical Engineering},

  volume={106},

  pages={777--784},

  year={2017},

  publisher={Elsevier}

}

@article{Barrault04,
author = {Barrault, Maxime and Maday, Yvon and Nguyen, Ngoc and Patera, Anthony},
year = {2004},
month = {11},
pages = {667–672},
title = {An ‘empirical interpolation’ method: Application to efficient reduced-basis discretization of partial differential equations},
volume = {339},
journal = {Comptes Rendus Mathematique},
doi = {10.1016/j.crma.2004.08.006}
}

@article{cstefuanescu2013pod,

  title={Nonlinear model order reduction of an ADI implicit shallow water equations model},

  author={{\c{S}}tef{\u{a}}nescu, R{\u{a}}zvan and Navon, Ionel Michael},

  journal={Journal of Computational Physics},

  volume={237},

  pages={95--114},

  year={2013},

  publisher={Elsevier}

}

@article{cicci2022deep,

  title={Deep HyROMnet: A deep learning-based operator approximation for hyper-reduction of nonlinear parametrized PDEs},

  author={Cicci, Ludovica and Fresca, Stefania and Manzoni, Andrea},

  journal={Journal of Scientific Computing},

  volume={93},

  number={2},

  pages={57},

  year={2022},

  publisher={Springer}

}

@article{cocola2023hyper,

  title={Hyper-Reduced Autoencoders for Efficient and Accurate Nonlinear Model Reductions},

  author={Cocola, Jorio and Tencer, John and Rizzi, Francesco and Parish, Eric and Blonigan, Patrick},

  journal={arXiv preprint arXiv:2303.09630},

  year={2023}

}

@article{argaud2018sensor,

  title={Sensor placement in nuclear reactors based on the generalized empirical interpolation method},

  author={Argaud, J-P and Bouriquet, Bertrand and De Caso, F and Gong, Helin and Maday, Yvon and Mula, Olga},

  journal={Journal of Computational Physics},

  volume={363},

  pages={354--370},

  year={2018},

  publisher={Elsevier}

}

@article{maday2008parareal,
  title={The parareal in time algorithm},
  author={Maday, Yvon},
  year={2008}
}

\end{document}